
\def\LaTeX{%
  \let\Begin\begin
  \let\End\end
  \let\salta\relax
  \let\finqui\relax
  \let\futuro\relax}

\def\UK{\def\our{our}\let\sz s}
\def\USA{\def\our{or}\let\sz z}

\UK



\LaTeX

\USA


\salta

\documentclass[twoside,12pt]{article}
\setlength{\textheight}{24cm}
\setlength{\textwidth}{16cm}
\setlength{\oddsidemargin}{2mm}
\setlength{\evensidemargin}{2mm}
\setlength{\topmargin}{-15mm}
\parskip2mm


\usepackage[usenames,dvipsnames]{color}
\usepackage{amsmath}
\usepackage{amsthm}
\usepackage{amssymb}
\usepackage[mathcal]{euscript}
\usepackage{cite}
%
%


\definecolor{viola}{rgb}{0.3,0,0.7}
\definecolor{ciclamino}{rgb}{0.5,0,0.5}

\def\pier #1{{\color{red}#1}}
\def\juerg #1{{\color{blue}#1}}
\def\gianni #1{{\color{ciclamino}#1}}
\def\juergen #1{{\color{red}#1}}

\def\pier #1{#1}
\def\juerg #1{#1}
\def\gianni #1{#1}
\def\juergen #1{#1}



\bibliographystyle{plain}


%

\finqui

\def\Beq{\Begin{equation}}
\def\Eeq{\End{equation}}
\def\Bsist{\Begin{eqnarray}}
\def\Esist{\End{eqnarray}}

\def\Bthm{\Begin{theorem}}
\def\Ethm{\End{theorem}}

\def\Bprop{\Begin{proposition}}
\def\Eprop{\End{proposition}}
\def\Bcor{\Begin{corollary}}
\def\Ecor{\End{corollary}}
\def\Brem{\Begin{remark}\rm}
\def\Erem{\End{remark}}

\def\Bnot{\Begin{notation}\rm}
\def\Enot{\End{notation}}

\def\Bcenter{\Begin{center}}
\def\Ecenter{\End{center}}
\let\non\nonumber




\def\step #1 \par{\medskip\noindent{\bf #1.}\quad}


\def\Lip{Lip\-schitz}
\def\Holder{H\"older}

\def\aand{\quad\hbox{and}\quad}

\def\lhs{left-hand side}
\def\rhs{right-hand side}
\def\sfw{straightforward}



\def\multibold #1{\def\arg{#1}%
  \ifx\arg\pto \let\next\relax
  \else
  \def\next{\expandafter
    \def\csname #1#1#1\endcsname{{\bf #1}}%
    \multibold}%
  \fi \next}

\def\pto{.}

\def\multical #1{\def\arg{#1}%
  \ifx\arg\pto \let\next\relax
  \else
  \def\next{\expandafter
    \def\csname cal#1\endcsname{{\cal #1}}%
    \multical}%
  \fi \next}


\def\multimathop #1 {\def\arg{#1}%
  \ifx\arg\pto \let\next\relax
  \else
  \def\next{\expandafter
    \def\csname #1\endcsname{\mathop{\rm #1}\nolimits}%
    \multimathop}%
  \fi \next}

\multibold
qwertyuiopasdfghjklzxcvbnmQWERTYUIOPASDFGHJKLZXCVBNM.

\multical
QWERTYUIOPASDFGHJKLZXCVBNM.

\multimathop
diag dist div dom mean meas sign supp .


\def\accorpa #1#2{\eqref{#1}--\eqref{#2}}
\def\Accorpa #1#2 #3 {\gdef #1{\eqref{#2}--\eqref{#3}}%
  \wlog{}\wlog{\string #1 -> #2 - #3}\wlog{}}


\def\separa{\noalign{\allowbreak}}

\def\somma #1#2#3{\sum_{#1=#2}^{#3}}
\def\tonde #1{\left(#1\right)}

\def\<#1>{\mathopen\langle #1\mathclose\rangle}
\def\norma #1{\mathopen \| #1\mathclose \|}
\def\Norma #1{\left\| #1 \right\|}

\def\[#1]{\mathopen\langle\!\langle #1\mathclose\rangle\!\rangle}

\def\iot {\int_0^t}
\def\ioT {\int_0^T}
\def\intQt{\int_{Q_t}}
\def\intQ{\int_Q}
\def\iO{\int_\Omega}

\def\dt{\partial_t}

\def\cpto{\,\cdot\,}

\def\checkmmode #1{\relax\ifmmode\hbox{#1}\else{#1}\fi}
\def\aeO{\checkmmode{a.e.\ in~$\Omega$}}
\def\aeQ{\checkmmode{a.e.\ in~$Q$}}

\def\aet{\checkmmode{a.e.\ in~$(0,T)$}}

\def\aat{\checkmmode{for a.a.~$t\in(0,T)$}}


\def\erre{{\mathbb{R}}}




\def\genspazio #1#2#3#4#5{#1^{#2}(#5,#4;#3)}
\def\spazio #1#2#3{\genspazio {#1}{#2}{#3}T0}

\def\L {\spazio L}
\def\H {\spazio H}
\def\W {\spazio W}

\def\C #1#2{C^{#1}([0,T];#2)}


\def\Lx #1{L^{#1}(\Omega)}
\def\Hx #1{H^{#1}(\Omega)}

\def\LQ #1{L^{#1}(Q)}

\def\Luno{\Lx 1}
\def\Ldue{\Lx 2}

\def\Huno{\Hx 1}
\def\Hdue{\Hx 2}
\def\Hunoz{{H^1_0(\Omega)}}


\def\LQ #1{L^{#1}(Q)}


\let\theta\vartheta
\let\eps\varepsilon
\let\phi\varphi

\let\hat\widehat

\let\TeXchi\chi                         
\newbox\chibox
\setbox0 \hbox{\mathsurround0pt $\TeXchi$}
\setbox\chibox \hbox{\raise\dp0 \box 0 }
\def\chi{\copy\chibox}


\def\pn{\mathbf{n}}

\def\VA #1{V_A^{#1}}
\def\VB #1{V_B^{#1}}
\def\VC #1{V_C^{#1}}

\def\Ar{{A,\,\rho}}
\def\Bs{{B,\,\sigma}}
\def\Ct{{C,\,\tau}}

\def\ful{f_1^\lambda}
\def\Ful{F_1^\lambda}
\def\Fl{F^\lambda}
\def\fl{f^\lambda}

\def\phil{\phi^\lambda}
\def\mul{\mu^\lambda}
\def\Sl{S^\lambda}

\def\phih{\hat\phi_h}
\def\muh{\hat\mu_h}
\def\Sh{\hat S_h}
\def\overphih{\overline\phi_h}
\def\overmuh{\overline\mu_h}
\def\overSh{\overline S_h}
\def\underphih{\underline\phi_h}
\def\undermuh{\underline\mu_h}
\def\underSh{\underline S_h}

\def\phin{\phi^n}
\def\phinp{\phi^{n+1}}

\def\mun{\mu^n}
\def\munp{\mu^{n+1}}

\def\Sn{S^n}
\def\Snp{S^{n+1}}
\def\dhphin{\frac{\phinp-\phin}h}
\def\dhmun{\frac{\munp-\mun}h}
\def\dhSn{\frac{\Snp-\Sn}h}

\def\phim{\phi^m}
\def\mum{\mu^m}
\def\Sm{S^m}

\def\phiz{\phi_0}
\def\muz{\mu_0}
\def\Sz{S_0}

\def\phiu{\phi^1}
\def\muu{\mu^1}
\def\Su{S^1}

\def\soluz{(\mu,\phi,S)}
\def\soluzl{(\mul,\phil,\Sl)}
\def\ds{\,ds}

\Begin{document}


%
\title{Well-posedness and regularity \\for a fractional tumor growth model}
\author{}
\date{}

\maketitle
\Bcenter
\vskip-1cm
{\large\sc Pierluigi Colli$^{(1)}$}\\
{\normalsize e-mail: {\tt pierluigi.colli@unipv.it}}\\[.25cm]
{\large\sc Gianni Gilardi$^{(1)}$}\\
{\normalsize e-mail: {\tt gianni.gilardi@unipv.it}}\\[.25cm]
{\large\sc J\"urgen Sprekels$^{(2)}$}\\
{\normalsize e-mail: {\tt sprekels@wias-berlin.de}}\\[.45cm]
$^{(1)}$
{\small Dipartimento di Matematica ``F. Casorati'', Universit\`a di Pavia}\\
{\small and Research Associate at the IMATI -- C.N.R. Pavia}\\
{\small via Ferrata 5, 27100 Pavia, Italy}\\[.2cm]
$^{(2)}$
{\small Department of Mathematics, Humboldt-Universit\"at zu Berlin}\\
{\small Unter den Linden 6, 10099 Berlin, Germany}\\
{\small and Weierstrass Institute for Applied Analysis and Stochastics}\\
{\small Mohrenstrasse 39, 10117 Berlin, Germany}
\Ecenter

%
\Begin{abstract}\noindent
In this paper, we study a system of three evolutionary operator equations involving fractional powers
of selfadjoint, monotone, unbounded, linear operators having compact resolvents. This system constitutes a
generalization of a phase field system of Cahn--Hilliard type modelling tumor growth that has been proposed
in Hawkins-Daarud et al. ({\em Int.\,\,J. Numer.\,\,Math.\,\,Biomed.\,\,Eng.\,\,{\bf 28} (2012), 3--24}) and  
investigated in recent papers co-authored by the present authors and E. Rocca. The model consists of a 
Cahn--Hilliard equation for the tumor cell fraction \,$\varphi$, coupled to a reaction-diffusion equation for a
function $\,S\,$ representing the nutrient-rich extracellular water volume fraction. Effects due to fluid
motion are neglected.  
The generalization investigated in this paper is motivated by the possibility that the diffusional regimes governing the evolution of the different constituents of the model may be of different (e.g., fractional) type. 
Under rather general assumptions, well-posedness and regularity results are shown. In particular, by 
writing the equation governing the evolution of the chemical potential in the form of a general variational 
inequality, also singular or nonsmooth contributions of
logarithmic or of double obstacle type to the energy density can be admitted.  
  
\End{abstract}
\salta
\pagestyle{myheadings}
\newcommand\testopari{\sc Colli \ --- \ Gilardi \ --- \ Sprekels}
\newcommand\testodispari{\sc Generalized fractional Cahn--Hilliard {}{tumor growth} system}
\markboth{\testopari}{\testodispari}
\finqui
%

\section{Introduction}
\label{Intro}
\setcounter{equation}{0}
Let $\Omega\subset\erre^3$ denote an open, bounded, and connected set
with smooth boundary $\Gamma$ and unit outward normal $\pn$, let $T>0$ be given,
and set $Q_t:=\Omega\times (0,t)$ for $t\in (0,T)$ and  $Q:=\Omega\times (0,T)$, as well as 
$\Sigma:=\Gamma\times (0,t)$. We investigate in this paper the evolutionary system
\Bsist
  && \alpha\, \dt\mu + \dt\phi + A^{2\rho} \mu 
  = P(\phi)(S-\mu){}{\qquad\mbox{in \,$Q$,}}
  \label{Iprima}
  \\
  && \mu =  \beta\, \dt\phi + B^{2\sigma}\phi + f(\phi) {}{\qquad\mbox{in \,$Q$}},
  \label{Iseconda}
  \\
  && \dt S + C^{2\tau}S 
  = - P(\phi)(S-\mu){}{\qquad\mbox{in \,$Q$,}}
  \label{Iterza}
  \\
  && \mu(0) = \muz, \quad
  \phi(0) = \phiz, \quad
  S(0) = \Sz,{}{\qquad\mbox{in \,$\Omega$.}}
  \label{Icauchy}
\Esist
In the above system, {}{$\alpha>0$ and $\beta>0$, and} $A^{2\rho}$, $B^{2\sigma}$, $C^{2\tau}$, with $r,\,\sigma,\,\tau>0$,
denote fractional powers of the selfadjoint,
monotone, and unbounded linear operators $A$, $B$, and~$C$, respectively, 
which are supposed to be densely defined in $H:=\Ldue$ and to have compact resolvents.
Moreover, $f$~denotes the derivative of a double-well potential~$F$.
Typical and physically significant examples of $F$ 
are the so-called {\em classical regular potential}, the {\em logarithmic potential\/},
and the {\em double obstacle potential\/}, which are given, in this order,~by
\begin{align}
\label{regpot}
  & F_{reg}(r) := \frac 14 \, (r^2-1)^2 \,,
  \quad r \in \erre,\\[1mm] 
  \label{logpot}
  &{}{ F_{log}(r) := \left\{\begin{array}{ll}
  \bigl( (1+r)\ln (1+r)+(1-r)\ln (1-r) \bigr) - c_1 r^2\,,&
  \quad r \in (-1,1)\\
  2\log(2)-c_1\,,&\quad r\in\{-1,1\}\\
  +\infty\,,&\quad r\not\in [-1,1]
\end{array}\right. \,,}
  \\[1mm]
\label{obspot}
  & F_{2obs}(r) :=  c_2\left(1 -  r^2\right)
  \quad \hbox{if $|r|\leq1$}
  \aand
  F_{2obs}(r) := +\infty
  \quad \hbox{if $|r|>1$}.
\end{align}
Here, the constants $c_i$ in \eqref{logpot} and \eqref{obspot} satisfy
$c_1>1$ and $c_2>0$, so that the corresponding functions are nonconvex.
In cases like \eqref{obspot}, one has to split $F$ into a nondifferentiable convex part~$F_1$ 
(the~indicator function of $[-1,1]$, in the present example) and a smooth perturbation~$F_2$.
Accordingly, in the term $f(\phi)$ appearing in~\eqref{Iseconda}, 
one has to replace the derivative $F_1'$ of the convex part $F_1$
by the subdifferential $f_1:=\partial F_1$ and interpret \eqref{Iseconda} as a differential inclusion
or as a variation inequality involving $F_1$ rather than~$f_1$.
Furthermore, the function $P$ {}{occurring} in \eqref{Iprima} and \eqref{Iterza}
is nonnegative and smooth.
Finally, the terms on the \rhs s in \eqref{Icauchy} are prescribed initial data.

The above system is a generalization of a phase field system of Cahn--Hilliard type~\cite{Cahn} modelling tumor growth. The original model \juergen{was} proposed in \cite{HZO12}, then extended 
in \cite{HKNZ, CGH15}, and investigated in \cite{CGRS1, CGRS2, FGR} from the 
viewpoint of well-posedness, regularity, \juergen{and} asymptotic analyses; 
instead, the papers \cite{CRW, CGRS3, Sig, S_a, S_b} \juergen{were} concerned with various
optimal control problems that have been set for this class of models. In the mentioned contributions, the three operators $A^{2\rho}$, $B^{2\sigma}$, $C^{2\tau}$ are nothing but the operator $ - \Delta$, with homogeneous Neumann boundary conditions. Concerning the meanings of the \juergen{variables} of the system \eqref{Iprima}--\eqref{Icauchy},  $\varphi$ represents an order parameter accounting for the tumor fraction,  
and $S$ stands for a nutrient concentration, while  
the third unknown $\mu$ is the related chemical potential, specified by \eqref{Iseconda} 
as for the viscous Cahn--Hilliard equation. Some interest of 
\eqref{Iprima}--\eqref{Iterza} 
becomes immediately evident and relies to the fact that we can admit different fractional 
operators in the description of the evolution of a tumor growth.  

Modelling the dynamics of tumor growth has recently become an important issue in applied mathematics (see, e.g.,~\cite{CL2010,WLFC}). Indeed, a noteworthy interest arose 
among mathematicians and applied scientists on the dynamics of tumor cells inside parts of the human body. Thus, a significant number of models have been introduced and discussed, with numerical simulations as well, in connection and comparison with 
the behavior of other special materials: one may see \cite{BLM,CLLW,CL2010,FBG2006,Fri2007,Lowen10,HDPZO,OHP,WLFC,WZZ}.

As diffuse interface models are concerned, we note that these models mostly 
use the Cahn--Hilliard framework, which is related to the theory of phase 
transitions, and which is used extensively in materials science and multiphase fluid flow.
Actually, one can distinguish between two main classes of models. 
The first one considers the tumor and healthy cells as inertialess fluids
including effects generated by fluid flow development, 
postulating a Darcy or a Brinkman law. To this concern, we refer to
\juergen{\cite{DFRSS,EGAR,FLRS,GL2016,GL2018,GLNS, GLSS, JWZ,LTZ,SW, Wise2011}} 
(see also \cite{BCG,ConGio,DG,FW2012,GioGrWu,WW2012, WZ2013} for local or nonlocal Cahn--Hilliard systems with Darcy or Brinkman law),
and we point out that further mechanisms such as chemotaxis 
and active transport can be taken into account.
The other class, \juergen{to} which the model leading to \eqref{Iprima}--\eqref{Icauchy} 
belongs, neglects the velocity and admits as variables concentrations and chemical potential.  Let us quote a group of contributions inside this class,
namely \cite{CRW,CWSL,CGMR, FLR, GL2017-1, GL2017-2,GLR, MRS}. 
To our knowledge, up to now fractional operators have not yet been 
dealt with in \juergen{either of these two} groups of models, although 
one may also wonder about nonlocal operators. 

All in all, fractional operators represent nowadays a challenging subject for 
mathematicians: they have been used in a number of situations, and 
there is already a wide literature about equations and systems with fractional terms.
In particular, different variants of fractional operators 
have been considered and employed. For a review of some related work, 
let us refer the interested reader to our recent papers 
\cite{CGS18, CGS19} and \cite{CG1}, which offer a recapitulation 
of various contributions. 
In our approach here, \juergen{which} follows closely the setting 
used in \cite{CG1,CGS18, CGS19, CGS22}, we deal with fractional 
operators defined via spectral theory. Then we can easily consider
powers of a second-order elliptic operator with either Dirichlet or Neumann or 
Robin \juergen{homogeneous} boundary conditions, as well as other operators like, e.g., 
fourth-order ones or systems involving the Stokes operator. 
The precise framework for our fractional operators
$A^{2\rho}$, $B^{2\sigma}$, $C^{2\tau}$, is given in the first part 
of Section~\ref{STATEMENT}.%

The remainder of the paper is organized as follows. 
In the next section, we list our assumptions and notations and state our results.
The uniqueness of the solution is proved in Section~\ref{UNIQUENESS},
while its existence is established {in Section~\ref{EXISTENCE}}.
The proof is prepared by the study of  
approximating discrete problems, which \juergen{are} introduced 
and solved in the subsections of the same section.
Finally, the last section is devoted to the regularity of the solution.


\section{Statement of the problem and results}
\label{STATEMENT}
\setcounter{equation}{0}

In this section, we state precise assumptions and notations and present our results.
As mentioned above, the set $\Omega\subset\erre^3$ is bounded, connected, and smooth,
with volume $|\Omega|$ and outward unit normal vector field $\pn$ on $\Gamma:=\partial\Omega$.
Moreover, $\partial_\pn$ stands for the corresponding normal derivative.
We set
\Beq
  H := \Ldue
  \label{defH}
\Eeq
and denote by $\norma\cpto$ and $(\cpto,\cpto)$ the standard norm and inner product of~$H$. Now, we start introducing our assumptions.
As for the operators, we first postulate that
\Bsist
  && A:D(A)\subset H\to H , \quad
  B:D(B)\subset H\to H
  \aand
  C:D(C)\subset H\to H
  \quad \hbox{are}
  \non
  \\
  && \hbox{unbounded, monotone, selfadjoint, linear operators with compact resolvents.} 
  \qquad
  \label{hpABC} 
\Esist
Therefore, there are sequences 
$\{\lambda_j\}$, $\{\lambda'_j\}$, $\{\lambda''_j\}$ and $\{e_j\}$, $\{e'_j\}$, $\{e''_j\}$ 
of eigenvalues and of corresponding eigenfunctions satisfying
\Bsist
  && A e_j = \lambda_j e_j, \quad
  B e'_j = \lambda'_j e'_j
  \aand
  C e''_j = \lambda''_j e''_j
  \non
  \\
  && \quad \hbox{with} \quad
  (e_i,e_j) = (e'_i,e'_j) = (e''_i,e''_j) = \delta_{ij}
  \quad \hbox{for $i,j=1,2,\dots$}
  \label{eigen}
  \\
  \separa
  && 0 \leq \lambda_1 \leq \lambda_2 \leq \dots , \quad
  0 \leq \lambda'_1 \leq \lambda'_2 \leq \dots
  \aand
  0 \leq \lambda''_1 \leq \lambda''_2 \leq \dots  
  \non
  \\
  && \quad \hbox{with} \quad
  \lim_{j\to\infty} \lambda_j
  = \lim_{j\to\infty} \lambda'_j
  = \lim_{j\to\infty} \lambda''_j
  = + \infty,
  \label{eigenvalues}
  \\[1mm]
  && \hbox{$\{e_j\}$, $\{e'_j\}$ and $\{e''_j\}$ are complete systems in $H$}.
  \label{complete}
\Esist
As a consequence, we can define the powers of these operators with arbitrary 
positive real exponents as done below.
As far as the first operator is concerned, we have for $r>0$
\Bsist
  && \VA\rho := D(A^\rho)
  = \Bigl\{ v\in H:\ \somma j1\infty |\lambda_j^\rho (v,e_j)|^2 < +\infty \Bigr\}
  \aand
  \label{defdomAr}
  \\[-3mm]
  && A^\rho v = \somma j1\infty \lambda_j^\rho (v,e_j) e_j
  \quad \hbox{for $v\in\VA\rho$},
  \label{defAr}
\Esist
the series being convergent in the strong topology of~$H$,
due to the properties \eqref{defdomAr} of the coefficients.
We endow $\VA\rho$ with the graph norm, i.e., we~set
\Beq
  (v,w)_\Ar := (v,w) + (A^\rho v,A^\rho w)
  \aand
  \norma v_\Ar := (v,v)_\Ar^{1/2}
  \quad \hbox{for $v,w\in\VA\rho$},
  \label{defnormaAr}
\Eeq
and obtain a Hilbert space.
In {}{the same} way, 
we can define the power $B^\sigma$ and $C^\tau$ for every $\sigma>0$ and $\tau>0$,
starting from \accorpa{hpABC}{complete} for $B$ and~$C$.
We therefore set
\Bsist
  && \VB\sigma := D(B^\sigma)
  \aand
  \VC\tau := D(C^\tau)
  \quad \hbox{with the norms $\norma\cpto_\Bs$ and $\norma\cpto_\Ct$}
  \non
  \\
  && \quad \hbox{associated with the inner products}
  \label{defBsCr}
  \non
  \\
  && (v,w)_\Bs := (v,w) + (B^\sigma v,B^\sigma w)
  \aand
  (v,w)_\Ct := (v,w) + (C^\tau v,C^\tau w),
  \non
  \\
  && \quad \hbox{for $v,w\in \VB\sigma$ and $v,w\in\VC\tau$, respectively}.
  \label{defprodBsCr}
\Esist
Since $\lambda_j\geq0$ for every~$j$, {}{one immediately deduces from the definition of $A^\rho$} that
\Bsist
  && A^\rho : \VA\rho \subset H \to H
  \quad \hbox{is maximal monotone and} 
  \non
  \\
  && {}{\eps I + A^\rho} : \VA\rho \to H
  \quad \hbox{is a {}{topological} isomorphism for every $\eps>0$},
  \label{maxmon}
\Esist
where $I:H\to H$ is the identity operator. Similar results hold 
for $B^\sigma$ and~$C^\tau$.
It is clear that, for every $\rho_1,\,\rho_2>0$, we have that
\Beq
  (A^{\rho_1+\rho_2} v,w)
  = (A^{\rho_1} v, A^{\rho_2} w)
  \quad \hbox{for every $v\in\VA{\rho_1+\rho_2}$ and $w\in\VA{\rho_2}$},
  \label{propA}
\Eeq
and that similar relations holds for the other two types of fractional operators.
Due to these properties, we can define proper extensions of the operators that
allow values in dual spaces.
In particular, we can write variational formulations of the equation \accorpa{Iprima}{Iterza}.
It is convenient to use the notations
\Beq
  \VA{-\rho} := (\VA\rho)^* , \quad
  \VB{-\sigma} := (\VB\sigma)^*, 
  \aand
  \VC{-\tau} := (\VC\tau)^* ,
  \quad \hbox{for $\pier{\rho},\,\sigma,\,\tau>0$}.
  \label{defnegspaces}
\Eeq
Thus, we have that
\Beq
  A^{2\rho} \in \calL(\VA\rho; \VA{-\rho}) , \quad
  B^{2\sigma} \in \calL(\VB\sigma; \VB{-\sigma}),
  \aand
  C^{2\tau} \in \calL(\VC\tau; \VC{-\tau}) ,
  \label{extensions2}
\Eeq
as well as
\Beq
  A^\rho \in \calL(H; \VA{-\rho}) , \quad
  B^\sigma \in \calL(H; \VB{-\sigma}),
  \aand
  C^\tau \in \calL(H; \VC{-\tau}) .
  \label{extensions1}
\Eeq
The symbol $\<\cpto,\cpto>_\Ar$ will be used for the duality pairing
between $\VA{-\rho}$ and~$\VA\rho$. 
Moreover, we identify $H$ with a subspace of $\VA{-\rho}$
in the usual way, i.e., such that
\Beq
  \< v,w >_\Ar = (v,w)
  \quad \hbox{for every $v\in H$ and $w\in\VA\rho$}.
  \label{identification}
\Eeq
Analogously, we have that $H\subset\VB{-\sigma}$ and $H\subset\VC{-\tau}$ and use similar notations.

From now on, we assume {}{that}
\Beq
  \hbox{\gianni{$\alpha$, $\beta$}, $\rho$, $\sigma$ and $\tau$ are fixed positive real numbers}.
  \label{hpexponents}
\Eeq
\gianni{Moreover, \pier{for} some of our results we have to require} 
the following continuous embeddings of Sobolev type:
\Beq
  \VA\rho \subset \Lx 4
  \aand
  \VC\tau \subset \Lx 4 \,.
  \label{embeddings}
\Eeq
\gianni{\juerg{Under} these assumptions}, we can choose some $M\geq1$ such that
\Beq
  \norma v_4 \leq M \, \norma v_\Ar
  \aand
  \norma v_4 \leq M \norma v_\Ct
  \label{normembeddings}
\Eeq
for every $v\in\VA\rho$ and $v\in\VC\tau$, respectively.

\Brem
\label{Embeddings}
For instance, the first embedding \eqref{embeddings} is satisfied 
if \gianni{$A=-\Delta$}, the \gianni{(negative)} Laplace operator, with domain $\Hdue\cap\Hunoz$ 
(thus, with homogeneous Dirichlet conditions, but similarly for the Neumann boundary conditions).
Indeed, \gianni{$\VA{1/2}=\Hunoz$} in this case.
Clearly, the same embedding holds true if $\rho$ is sufficiently close to~\gianni{$1/2$}.
\Erem

For the  nonlinear functions entering the equations \accorpa{Iprima}{Iterza} of our system, we postulate the properties listed below.
The notation $\Ful$ stands for the Moreau--Yosida regularization of~$F_1$ at the level~$\lambda>0$
(see, e.g., \cite[p.~39]{Brezis}).
\Bsist
  && F := F_1 + F_2, \quad \hbox{where:}
  \label{hpF}
  \\
  && F_1 : \erre \to [0,+\infty]
  \quad \hbox{is convex, proper, and l.s.c.,\ with} \quad
  F_1(0) = 0\pier{;}
  \label{hpFuno}
  \\
  \separa
  && F_2 : \erre \to \erre
  \quad \hbox{is of class $C^1$ with a \Lip\ continuous first derivative\pier{;}
}
  \qquad
  \label{hpFdue}
  \\
  && \Ful(s) + F_2(s) \geq - C_0
  \quad \hbox{for some constant $C_0$ and every $s\in\erre$\pier{;}
}
  \label{hpbelow}
  \\
  && P :\erre \to [0,+\infty) \quad \hbox{is bounded and \Lip\ continuous}.
  \label{hpP}
\Esist
\Accorpa\Hpnonlin hpF hpP

\Brem
\label{Hpbelow}
The assumption \eqref{hpbelow} can be supposed to hold just for sufficiently small $\lambda>0$.
A sufficient condition for this 
(see \cite[formula (3.1)]{CGS22} for some explanation)
is that $F$ satisfies an inequality of type
\Beq
  F(s) \geq c_1 s^2 - c_2\,,
  \quad \hbox{for some constants $c_i>0$ and every $s\in\erre$}.
  \label{pier}
\Eeq
Hence, \accorpa{hpFuno}{hpbelow} are 
fulfilled by all of the important potentials \accorpa{regpot}{obspot}.
\Erem

We set, for convenience,
\Beq
  f_1 := \partial F_1 
  \aand
  f_2 := F_2' \,.
  \label{deffunodue}  
\Eeq
\Accorpa\Hpstruttura hpexponents deffunodue
Moreover, we term $D(F_1)$ and $D(f_1)$ the effective domains of $F_1$ and~$f_1$, respectively.
We notice that $f_1$ is a maximal monotone graph in $\erre\times\erre$
and use the same symbol $f_1$ for the maximal monotone operators induced in $L^2$ spaces.
Observe that $D(F_1)=D(f_1)=\erre$ for $F=F_{reg}$, while $D(F_1)=[-1,1]$ and $D(f_1)=(-1,1)$
for \gianni{$F=F_{log}$. Finally, \pier{we have that} $D(F_1)=D(f_1)=[-1,1]$ if $F=F_{2obs}$}.

On account of \eqref{propA} and its analogues for~$B$ and~$C$, we give a weak formulation of the equations \accorpa{Iprima}{Iterza}.
Moreover, we present \eqref{Iseconda} as a variational inequality.
For the data, we make the following assumptions:
\Bsist
  && \muz \in H , \quad
  \phiz \in \VB\sigma 
  \quad \hbox{with} \quad
  F_1(\phiz) \in \Luno , 
  \aand
  \Sz \in H .
  \label{hpdati}
\Esist
{}{We then} look for a triplet $\soluz$ satisfying
\Bsist
  && \mu \in \H1{\VA{-\rho}} \cap \L\infty H \cap \L2{\VA\rho},
  \label{regmu}
  \\
  && \phi \in \H1H \cap \L\infty{\VB\sigma},
  \label{regphi}
  \\
  && S \in \H1{\VC{-\tau}} \cap \L\infty H \cap \L2{\VC\tau},
  \label{regS}
  \\
  && F_1(\phi) \in \LQ1,
  \label{regFuphi}
\Esist
\Accorpa\Regsoluz regmu regFuphi
and solving the system
\Bsist
  && \alpha \< \dt\mu(t) , v >_\Ar
  + \bigl( \dt\phi(t) , v \bigr)
  + ( A^\rho \mu(t) , A^\rho v )
  = \bigl( \gianni{P(\phi(t))} ( S(t) - \mu(t) ) , v \bigr)
  \non
  \\
  && \quad \hbox{for every $v\in\VA\rho$ and \aat,}
  \label{prima}
  \\[2mm]
  \separa
  && \beta \bigl( \dt\phi(t) , \phi(t) - v \bigr)
  + \bigl( B^\sigma \phi(t) , B^\sigma( \phi(t)-v) \bigr)
  \non
  \\
  && \quad {}
  + \iO F_1(\phi(t))
  + \bigl( f_2(\phi(t)) ,  \phi(t)-v \bigr)
  \leq \bigl( \mu(t) , \phi(t)-v \bigr)
  + \iO F_1(v)
  \non
  \\
  && \quad \hbox{for every $v\in\VB\sigma$ and \aat,}
  \label{seconda}
  \\[2mm]
  \separa
  && \< \dt S(t) , v >_\Ct
  + ( C^\tau S(t) ,C^\tau v )
  = - \bigl( \gianni{P(\phi(t))} ( S(t) - \mu(t) ) , v \bigr)
  \non
  \\
  && \quad \hbox{for every $v\in\VC\tau$ and \aat},
  \label{terza}
  \\
  && \mu(0) = \muz \,, \quad
  \phi(0) = \phiz\,, 
  \aand
  S(0) = \Sz \,.
  \label{cauchy}
\Esist
\Accorpa\Pbl prima cauchy
{}{Here, it is understood that $\,\,\iO F_1(v)=+\infty\,$ whenever $\,F_1(v)\not\in\gianni\Luno$.}

We notice at once that \eqref{seconda} is equivalent to its time-integrated variant, that is,
\Bsist
  && \beta \ioT \bigl( \dt\phi(t) , \phi(t) - v(t) \bigr) \, dt 
  + \ioT \bigl( B^\sigma\phi(t) , B^\sigma(\phi(t)-v(t)) \bigr)\,dt 
  \non
  \\ 
  && \quad {}
  + \intQ F_1(\phi) 
  + \ioT \bigl( f_2(\phi(t)) , \phi(t)-v(t) \bigr) \, dt 
  \non
  \\
  && \leq \ioT \bigl( \mu(t) , \phi(t)-v(t) \bigr)\,dt
  + \intQ F_1(v)
  \quad \hbox{for every $v\in\L2{\VB\sigma}$}.
  \qquad
  \label{intseconda}
\Esist
Here is our well-posedness and continuous dependence result.

\Bthm
\label{Wellposedness}
Let the assumptions \gianni{\eqref{hpABC}, \eqref{hpexponents}, \Hpnonlin\ and \eqref{deffunodue}}
on the structure of the system,
and \eqref{hpdati} on the data be fulfilled.
Then there exists at least one triplet $\soluz$ 
satisfying \Regsoluz\ and solving problem \Pbl. 
Moreover, for this solution we have the estimates
\Bsist
  && \alpha^{1/2} \norma\mu_{\L\infty H}
  + \norma{A^\rho\mu}_{\L2H}
  \non
  \\
  && \quad {}
  + \beta^{1/2} \norma{\dt\phi}_{\L2H}
  + \norma{B^\sigma\phi}_{\L\infty H}
  + \pier{\norma{F(\phi)}_{\L\infty{L^1(\Omega)}}}
  \non
  \\[1mm]
  && \quad {}
  + \norma S_{\L\infty H}
  + \norma{C^\tau S}_{\L2H}
  + \norma{P^{1/2}(\phi)(S-\mu)}_{\L2H}
  \non
  \\[1mm]
  && \juerg{\leq C_1 \bigl(
    \alpha^{1/2} \norma\muz
    + \norma{B^\sigma\phiz}
    + \norma{F(\phiz)}_{\Luno}
    + \norma\Sz
    + 1
  \bigr)},  \label{stimasoluz}
  \\[2mm]
  && \gianni{\norma{\dt(\alpha\mu+\phi)}_{\L2{\VA{-\rho}}}
  + \norma{\dt S}_{\L2{\VC{-\tau}}}}
  \non
  \\
  && \gianni{{}\leq C_2 \bigl(
    \alpha^{1/2} \norma\muz
    + \norma{B^\sigma\phiz}
    + \norma{F(\phiz)}_{\Luno}
    + \norma\Sz
    + 1
  \bigr)},
  \label{stimadtsoluz}
\Esist
\gianni{with a constant $C_1$ that depends only on~$\Omega$ and the constant $C_0$ from~\eqref{hpbelow},
and a constant $C_2$ that also depends on~$P$.
If, in addition, \eqref{pier} is satisfied, then we also have
\Beq
  \norma\phi_{\L\infty{\VB\sigma}}
  \leq \gianni{C_3} \bigl(
    \alpha^{1/2} \norma\muz
    + \norma{B^\sigma\phiz}
    + \norma{F(\phiz)}_{\Luno}
    + \norma\Sz
    + 1
  \bigr),
  \label{dapier}
\Eeq
where \gianni{$C_3$ depends on $\Omega$ and the constants $C_0$, $c_1$ and $c_2$ from \eqref{hpbelow} and~\eqref{pier}}.
Finally, the solution $\soluz$ is unique if the spaces $\VA\rho$ and $\VC\tau$ satisfy~\eqref{embeddings}.}
\Ethm

\Brem
\label{Controls}
More generally, we could add known forcing terms $u_\mu$, $u_\phi$ and $u_S$
to the \rhs s of equations \eqref{Iprima}, \eqref{Iseconda} and~\eqref{Iterza}, respectively,
and accordingly modify the definition of solution.
If we assume that
\Beq
  u_\mu ,\ u_\phi \,, u_S \in \L2H\,,
  \label{hpcontrols}
\Eeq
then we have a similar well-posedness result.
In estimate~\eqref{stimasoluz}, one has to modify the \rhs\ by adding 
the norms corresponding to \eqref{hpcontrols} 
(possibly multiplied by negative powers of $\alpha$ and~$\beta$).
This remark is useful if one has in mind to perform a control theory on the above system
with distributed controls.
\Erem

\gianni{Our next aim is \juerg{to prove} further \pier{properties} of the solution.
Indeed, from one side, one wishes to improve the regularity requirements \Regsoluz\
under suitable new assumptions on the data.
On the other hand, one wishes to have an equation
(or~at least a differential inclusion)
in place of the variational inequality~\eqref{seconda}.
The next results deal with these problems independently from each other, in principle.
The first one requires just something more on the initial data, indeed.
Namely, we assume that
\Beq
  \muz \in \VA\rho \,, \quad
  \phiz \in \VB{2\sigma}
  \quad \hbox{with} \quad
  f_1^\circ(\phiz) \in H ,
  \aand
  \Sz \in \VC\tau\,, 
  \label{hppiuregdati}
\Eeq
where for $s\in D(f_1)$,
the symbol $f_1^\circ(s)$ stands for the element of $f_1(s)$ having minimum modulus.
We notice that \eqref{hppiuregdati} implies \eqref{hpdati}
since $F_1(s)\leq F_1(0)+s f_1^\circ(s)$ for every $s\in D(f_1)$ by convexity and $F_1(0)=0$.
Hence the existence of a solution is still ensured.}

\gianni{For the other problem, one cannot expect anything that is similar to \eqref{Iseconda},
since no estimate for $f_1(\phi)$ is available in the general case.
\juerg{However}, if the assumptions on the structure are reinforced, \juerg{then}
one can recover~\eqref{Iseconda} at least as a differential inclusion.
The crucial condition is the following:
\begin{align}
  & \psi(v) \in H
  \aand
  \bigl( B^{2\sigma} v , \psi(v) \bigr) \geq 0,
  \quad \hbox{for every $v\in\VB{2\sigma}$}
  \,\hbox{ and every monotone}
  \non
  \\
  & \mbox{\pier{and} \Lip\ continuous function $\psi:\erre\to\erre$ \gianni{vanishing at the origin}}.
   \label{hpsecondeq}
\end{align}
We notice that this assumption is fulfilled if $B^{2\sigma}=-\Delta$ 
with zero Neumann boundary conditions.
Indeed, in this case it results that $\VB{2\sigma}=\{v\in\Hdue:\ \partial_\pn v=0\}$,
and, for every $\psi$ as in \eqref{hpsecondeq} and $v\in\VB{2\sigma}$, 
we have that $\psi(v)\in\Huno$ (since $v\in\Huno$)~\juerg{as well~as}
\Beq
  \bigl( B^{2\sigma} v, \psi(v) \bigr)
  = \iO (-\Delta v) \, \psi(v)
  = \iO \nabla v \cdot \nabla\psi(v)
  = \iO \psi'(v) |\nabla v|^2
  \geq 0 .
  \non
\Eeq
More generally, in place of the Laplace operator we can take the principal part
of an elliptic operator in divergence form with \Lip\ continuous coefficients,
provided that the normal derivative is replaced by the conormal derivative.
In any case, we can take the (zero) Dirichlet boundary conditions instead of the Neumann boundary conditions,
since the functions $\psi$ for which \eqref{hpsecondeq} is required \pier{has to satisfy $\psi(0)=0$.}}

\gianni{%
\Bthm
\label{Regularity}
Let the assumptions \eqref{hpABC}, \accorpa{hpexponents}{hpP}, and \eqref{deffunodue}
on the structure of the system be fulfilled.
Moreover, let the data satisfy~\eqref{hppiuregdati}.
Then the unique solution $\soluz$ to problem \Pbl\ enjoys the further regularity
\Bsist
  && \mu \in \H1H \cap \L\infty{\VA\rho} \cap \L2{\VA{2\rho}},
  \label{piuregmu}
  \\
  && \phi \in \W{1,\infty}H \cap \H1{\VB\sigma},
  \label{piuregphi}
  \\
  && S \in \H1H \cap \L\infty{\VC\tau} \cap \L2{\VC{2\tau}} .
  \label{piuregS}
\Esist
\Ethm
}%

\Bthm
\label{Secondeq}
Besides the assumptions \eqref{hpABC}, \eqref{hpexponents}, \Hpnonlin,  and \eqref{deffunodue}
on the structure of the system, let \eqref{hpsecondeq} be fulfilled.
Moreover, let the data satisfy~\eqref{hpdati}.
Then there exist a solution $\soluz$ to problem \Pbl\ and some $\xi$ such~that
\Bsist
  && \phi \in \L2{\VB{2\sigma}}
  \aand
  \xi \in \L2H\,,
  \label{regseconda}
  \\[1mm]
  && \beta \, \dt\phi + B^{2\sigma}\phi + \xi + f_2(\phi) = \mu
  \aand
  \xi \in f_1(\phi)
  \quad \aeQ \,.
  \label{eqseconda}
\Esist
Furthermore, even $\xi$ is unique under the further assumption \eqref{embeddings}.
\Ethm

\gianni{%
\Bcor
\label{Fullregularity}
Assume \eqref{hpABC}, \accorpa{hpexponents}{hpP}, \eqref{deffunodue}, and~\eqref{hpsecondeq}
for the structure of the system, 
and \eqref{hppiuregdati} for the data.
Then the unique solution $\soluz$ to problem \Pbl\ and the corresponding $\xi$ 
satisfy \accorpa{piuregmu}{eqseconda} as well~as
\Beq
  \phi \in \L\infty{\VB{2\sigma}}
  \aand
  \xi \in \L\infty H.
  \label{fullregularity}
\Eeq
\Ecor
}%

In the following, we make use of
the elementary identity and of the Young inequality
\Bsist
  \hskip-1cm && a (a-b)
  = \frac 12 \, a^2
  + \frac 12 \, (a-b)^2
  - \frac 12 \, b^2
  \quad \hbox{for every $a,b\in\erre$},
  \label{elementare}
  \\
  \hskip-1cm && ab\leq \delta a^2 + \frac 1 {4\delta}\,b^2
  \quad \hbox{for every $a,b\in\erre$ and $\delta>0$}.
  \label{young}
\Esist
\Accorpa\Elementari elementare young
Moreover, if $V$ is a Banach space and $v$ is any function in~$\L2V$, then 
we define $1*v\in\H1V$ by setting
\Beq
  (1*v)(t) := \iot v(s) \ds
  \quad \hbox{for $t\in[0,T]$}.
  \label{convoluz}
\Eeq
Notice that, for every $t\in[0,T]$, we have
\Bsist
  && \norma{(1*v)(t)}_V^2
  \leq T \iot \norma{v(s)}_V^2 \ds\,,
  \label{disconvoluz}
  \\
  && \bigl( 1*(uv) \bigr)(t)
  = u(t) (1*v)(t)
  - \iot u'(s) (1*v)(s) \ds
  \quad \hbox{if \ $u\in H^1(0,T)$}.
  \qquad
  \label{idconvoluz}
\Esist
As for the notation concerning norms, we use the symbol $\norma\cpto_V$
for the norm in the generic Banach space~$V$ (as~done in~\eqref{disconvoluz})
with the following exceptions:
the simpler symbol $\norma\cpto$ denotes the norms in~$H$, as already said;
for the norms in the spaces $\VA\rho$, $\VB\sigma$ and $\VC\tau$ 
we use the notations introduced above;
if \pier{$1\leq q\leq\infty$}, the norm in any $L^q$ space is denoted 
by~$\norma\cpto_q$.

Finally, we state a general rule that we follow throughout the paper 
as far as the constants are concerned.
We use a small-case italic $c$ without subscripts
for~different constants that may only depend on the final time~$T$,
the operators $A^\rho$, $B^\sigma$ and~$C^\tau$,  
the shape of the nonlinearities $F$ and~$P$,
and the properties of the data involved in the statements at hand.
The values of such constants might change from line to line 
and even within the same formula or chain of inequalities.
The symbol $c_\delta$ stand for (possibly different) constants
that depend on the parameter $\delta$ in addition.
It is clear that $c$ and $c_\delta$ do not depend 
on the regularization parameter $\lambda$ and the time step $h$ we introduce in the next sections.
With the aim of performing some asymptotic analyses 
as the parameters $\alpha$ and/or $\beta$ tend to zero,
we clearly specify that the values of $c$ or $c_\delta$
do not depend on $\alpha$ or~$\beta$.
Constants (possibly different from each other) that depend, e.g., on both $\alpha$ and $\beta$ are denoted by 
$c_{\alpha,\beta}$.
In contrast, we use different symbols 
(like $M$ in \eqref{normembeddings} or $C_0$ in~\eqref{hpbelow})
for precise values of constants that we want to refer~to.


\section{Uniqueness}
\label{UNIQUENESS}
\setcounter{equation}{0}

In this section, we give the proof of the uniqueness part of Theorem~\ref{Wellposedness}.
We pick two solutions $(\mu_i,\phi_i,S_i)$, $i=1,2$, 
and set for convenience
$\mu:=\mu_1-\mu_2$, $\phi:=\phi_1-\phi_2$ and $S:=S_1-S_2$.
Now, we write equation \eqref{prima} for these solutions 
and integrate the difference with respect to time.
Then, we test the equality {}{thus obtained by}~$\mu$.
At the same time, we test the difference of \eqref{terza} written for the two solutions by~$S$.
Then, we sum up, integrate over $(0,t)$ with an arbitrary $t\in(0,T]$, and rearrange.
The {}{resulting \lhs\ contains the term} \,\,$\frac12\norma{A^\rho(1*\mu)(t)}^2$.
Thus, we add the same quantity
\Beq
  \frac 12 \, \norma{(1*\mu)(t)}^2
  = \iot \bigl( (1*\mu)(s) , \mu(s) \bigr) \ds
  \non
\Eeq
to both sides, by choosing the former expression for the \lhs\ and the latter for the \rhs,
and use the definition \eqref{defnormaAr} with $v=(1*\mu)(t)$.
Similarly, we add the same quantity to both sides
in order to reconstruct the full norm~$\norma{S(\cpto)}_\Ct$ in the corresponding integral.
We \juerg{then} obtain the identity
\Bsist
  && \alpha \iot \norma{\mu(s)}^2 \ds
  + \iot \bigl( \phi(s) , \mu(s) \bigr) \ds
  + \frac 12 \, \norma{(1*\mu)(t)}_\Ar^2 
  \non
  \\
  && \quad {}
  + \frac 12 \, \norma{S(t)}^2
  + \iot \norma{S(s)}_\Ct^2 \ds
  \non
  \\
  \separa
  && = \iot \Bigl(
    \bigl(
      1 * [ P(\phi_1) (S_1 - \mu_1) - P(\phi_2) (S_2 - \mu_2) ]
    \bigr) (s) , \mu(s) 
  \Bigr) \ds
  \non
  \\
  \separa
  && \quad {}
  - \iot \Bigl(
    \bigl(
      P(\phi_1) (S_1 - \mu_1) - P(\phi_2) (S_2 - \mu_2)
    \bigr) (s) , S(s) 
  \Bigr) \ds 
  \non
  \\
  \separa
  && \quad {}
  + \iot \bigl( (1*\mu)(s) , \mu(s) \bigr) \ds
  + \iot \norma{S(s)}^2 \ds \,.
  \label{diff13}
\Esist
Now, we treat the first term of the \rhs.
In the sequel, $\delta$ is a positive parameter.
By an integration by parts, we~get that
\Bsist
  && \iot \Bigl(
    \bigl(
      1 * [ P(\phi_1) (S_1 - \mu_1) - P(\phi_2) (S_2 - \mu_2) ]
    \bigr) (s) , \mu(s) 
  \Bigr) \ds
  \non
  \\
  \separa
  && = - \iot \Bigl(
    \bigl(
     P(\phi_1) (S_1 - \mu_1) - P(\phi_2) (S_2 - \mu_2)
    \bigr) (s) , (1*\mu)(s) 
  \Bigr) \ds
  \non
  \\
  && \quad {}
  + \Bigl(
    \bigl(
      1 * [ P(\phi_1) (S_1 - \mu_1) - P(\phi_2) (S_2 - \mu_2) ]
    \bigr) (t) , (1*\mu)(t) 
  \Bigr)\,,
  \non
\Esist
and we denote by $Y_1$ and $Y_2$, in this order, the summands on the \rhs.
To handle these terms, we owe to the \Holder\ and Young inequalities,
the boundedness and the \Lip\ continuity of~$P$, and the embeddings~\eqref{embeddings}.
As for~$Y_1$, we~have
\Bsist
  && Y_1 
  = - \iot \bigl(
    [ P(\phi_1) (S-\mu) ](s) , (1*\mu)(s)
  \bigr) \ds
  \non
  \\
  && \qquad {}
  - \iot \bigl(
    [ (P(\phi_1)-P(\phi_2)) (S_2-\mu_2) ](s) , (1*\mu)(s)
  \bigr) \ds
  \non
  \\
  \separa
  && \quad \,\leq \,\delta \iot \bigl( \norma{S(s)}^2 + \norma{\mu(s)}^2 \bigr) \ds
  + c_\delta \iot \norma{(1*\mu)(s)}^2 \ds
  \non
  \\
  && \qquad {}
  + c \iot \norma{\phi(s)}_2 \bigl( \norma{S_2(s)}_4 + \norma{\mu_2(s)}_4 \bigr) \norma{(1*\mu)(s)}_4 \ds
  \non
  \\
  \separa
  && \quad\,\leq \,\delta \iot \bigl( \norma{S(s)}^2 + \norma{\mu(s)}^2 \bigr) \ds
  + c_\delta \iot \norma{(1*\mu)(s)}^2 \ds
  \non
  \\
  && \qquad {}
  + \delta \iot \norma{\phi(s)}_\Bs^2 \ds
  + c_\delta \iot \bigl( \norma{S_2(s)}_\Ct^2 + \norma{\mu_2(s)}_\Ar^2 \bigr) \norma{(1*\mu)(s)}_\Ar^2 \ds\,,
  \qquad
  \label{stimaY1}
\Esist
and we notice that the function $s\mapsto\norma{S_2(s)}_\Ct^2 + \norma{\mu_2(s)}_\Ar^2$
belongs to $L^1(0,T)$,  thanks to the regularity assumed for the solution $(\mu_2,\phi_2,S_2)$.
In order to deal with~$Y_2$,
we prepare  an estimate of a delicate term with the help of~\eqref{idconvoluz}.
Since $P$ is nonnegative, one of the {}{resulting terms turns out to be} nonpositive.
Thus, on account of \Holder's inequality, \eqref{disconvoluz}, and~\eqref{embeddings}, we~have
\Bsist
  && - \bigl(
    (1*[ P(\phi_1) \mu ])(t) , (1*\mu)(t)
  \bigr)
  \non
  \\
  && = \,- \bigl( P(\phi_1(t)) (1*\mu)(t) , (1*\mu)(t) \bigr)
  + \Bigl( \textstyle\iot P'(\phi_1(s)) \dt\phi_1(s)(1*\mu)(s) \ds , (1*\mu)(t) \Bigr)
  \non
  \\
  \separa
  && \leq\, \delta\, \norma{(1*\mu)(t)}_4^2 
  + c_\delta \Norma{\textstyle\iot P'(\phi_1(s)) \dt\phi_1(s) (1*\mu)(s) \ds}_{4/3}^2
  \non
  \\
  && \leq \,\delta\, \norma{(1*\mu)(t)}_4^2 
  + c_\delta \iot \norma{\dt\phi_1(s)}_2^2 \, \norma{(1*\mu)(s)}_4^2 \ds
  \non
  \\
  \separa
  && \leq \,\delta\, M \, \norma{(1*\mu)(t)}_\Ar^2 
  + c_\delta \iot \norma{\dt\phi_1(s)}^2 \, \norma{(1*\mu)(s)}_\Ar^2 \ds\,,
  \non
\Esist
where we observe that the function $s\mapsto\norma{\dt\phi_1(s)}^2\,$ belongs to $L^1(0,T)$.
At this point, we can estimate the term~$Y_2$ 
by using \eqref{disconvoluz}, this inequality, and \eqref{normembeddings} \pier{with $M\geq 1$}.
We have
\begin{align}
   Y_2\,
  &=\, \bigl( 
    ( 1*[P(\phi_1) (S-\mu) - (P(\phi_1)-P(\phi_2)) (S_2-\mu_2)])(t)
    , (1*\mu)(t)
  \bigr)
  \non
  \\[2mm]
  & \leq\, \delta \, \norma{(1*\mu)(t)}^2
  + c_\delta \,\norma{(1*[ P(\phi_1) S ])(t)}^2
  \non
  \\[2mm]
  & \qquad {}
  - \bigl( (1*[ P(\phi_1) \mu ])(t) , (1*\mu)(t)
  \bigr)
  \non
  \\[2mm]
  & \qquad {}
  + \delta \, \norma{(1*\mu)(t)}_4^2
  + c_\delta \norma{(1*[ (P(\phi_1)-P(\phi_2)) (S_2-\mu_2) ])(t)}_{4/3}^2 
  \non
  \\
 \separa
  &\, \leq\, \delta \, \norma{(1*\mu)(t)}^2
  + c_\delta \iot \norma{(P(\phi_1) S)(s)}^2 \ds
  \non
  \\ 
  & \qquad {}
  + \delta \,M \, \norma{(1*\mu)(t)}_\Ar^2 
  + c_\delta \iot \norma{\dt\phi_1(s)}^2 \, \norma{(1*\mu)(s)}_\Ar^2 \ds
  \non
  \\
  & \qquad {}
  + \delta \, \norma{(1*\mu)(t)}_4^2
  + c_\delta \iot \norma{[ (P(\phi_1)-P(\phi_2)) (S_2-\mu_2) ](s)}_{4/3}^2 \ds
  \non
  \\
  & \,\leq\, 3\,\delta\, M \, \norma{(1*\mu)(t)}_\Ar^2 
  + c_\delta \iot \norma{S(s)}^2 \ds
  + c_\delta \iot \norma{\dt\phi_1(s)}^2 \, \norma{(1*\mu)(s)}_\Ar^2 \ds
  \non
  \\
  &\qquad {}
  + c_\delta \iot \bigl( \norma{S_2(s)}_4^2 + \norma{\mu_2(s)}_4^2 \bigr) \norma{\phi(s)}^2 \ds\,,
  \label{stimaY2}
\end{align}
where the function $s\mapsto \norma{S_2(s)}_4^2 + \norma{\mu_2(s)}_4^2\,$
is known to belong to $L^1(0,T)$.
Indeed, we have \,$S_2\in\L2{\VC\tau}\subset\L2{\Lx4}$\, and \,$\mu_2\in\L2{\VA\rho}\subset\L2{\Lx4}$.

Now, we come back to \eqref{diff13} and estimate the second term on the \rhs,
which we call $Y_3$ for simplicity.
We~have
\begin{align}
   Y_3\,&
  = \,
  - \iot \Bigl(
    \bigl(
      P(\phi_1) (S - \mu)
      + (P(\phi_1) - P(\phi_2)) (S_2 - \mu_2)
    \bigr) (s) , S(s) 
  \Bigr) \ds 
  \non
  \\
  \separa
  & \leq\, \delta \iot \norma{\mu(s)}^2 \ds
  + c_\delta \iot \norma{S(s)}^2 \ds
  \non
  \\
  &\qquad {}
  + \delta \iot \norma{S(s)}_4^2 \ds
  + c_\delta \iot \bigl( \norma{S_2(s)}_4^2 + \norma{\mu_2(s)}_4^2 \bigr) \norma{\phi(s)}^2 \ds
  \non
  \\
  & \leq \,\delta \iot \norma{\mu(s)}^2 \ds
  + c_\delta \iot \norma{S(s)}^2 \ds
  \non
  \\
  & \qquad {}
  + \delta M \iot \norma{S(s)}_\Ct^2 \ds
  + c_\delta \iot \bigl( \norma{S_2(s)}_\Ct^2 + \norma{\mu_2(s)}_\Ar^2 \bigr) \norma{\phi(s)}^2 \ds\,.
  \label{stimaY3}
\end{align}
At this point, we recall \accorpa{diff13}{stimaY3} and
use the Schwarz and Young inequalities to estimate the first term of the last line of~\eqref{diff13},  
{}{in order to get the estimate}
\begin{align}
  & \alpha \iot \norma{\mu(s)}^2 \ds
  + \iot \bigl( \phi(s) , \mu(s) \bigr) \ds
  + \frac 12 \, \norma{(1*\mu)(t)}_\Ar^2 
  \non
  \\
  & \quad {}
  + \frac 12 \, \norma{S(t)}^2
  + \iot \norma{S(s)}_\Ct^2 \ds
  \non
  \\
  \separa
  & \leq\, \delta \iot \bigl( \norma{S(s)}^2 + \norma{\mu(s)}^2 \bigr) \ds
  + c_\delta \iot \norma{(1*\mu)(s)}^2 \ds
  \non
  \\
  & \quad {}
  + \delta \iot \norma{\phi(s)}_\Bs^2 \ds
  + c_\delta \iot \bigl( \norma{S_2(s)}_\Ct^2 + \norma{\mu_2(s)}_\Ar^2 \bigr) \norma{(1*\mu)(s)}_\Ar^2 \ds
  \non
  \\
  & \quad {}
  + 3\,\delta\, M \, \norma{(1*\mu)(t)}_\Ar^2 
  + c_\delta \iot \norma{S(s)}^2 \ds
  + c_\delta \iot \norma{\dt\phi_1(s)}^2 \, \norma{(1*\mu)(s)}_\Ar^2 \ds
  \non
  \\
  & \quad {}
  + c_\delta \iot \bigl( \norma{S_2(s)}_4^2 + \norma{\mu_2(s)}_4^2 \bigr) \norma{\phi(s)}^2 \ds
  + \delta \iot \norma{\mu(s)}^2 \ds
  + c_\delta \iot \norma{S(s)}^2 \ds
  \non
  \\
  \separa
  & \quad {}
  + \delta M \iot \norma{S(s)}_\Ct^2 \ds
  + c_\delta \iot \bigl( \norma{S_2(s)}_\Ct^2 + \norma{\mu_2(s)}_\Ar^2 \bigr) \norma{\phi(s)}^2 \ds
  \non
  \\
  & \quad {}
  + \delta \iot \norma{\mu(s)}^2 \ds
  + c_\delta \iot \norma{(1*\mu)(s)}^2 \ds
  + \iot \norma{S(s)}^2 \ds \,.
  \label{dadiff13}
\end{align}
Next, we use the variational inequality~\eqref{seconda},
writing it for the two solutions, and testing the resulting inequalities by~$\phi_2$ and $\phi_1$, respectively.
Now, we sum up and notice that the contributions involving $F_1$ cancel out.
By integrating over~$(0,t)$, using the \Lip\ continuity of~$f_2$, 
and adding the same term to both sides in order to recover the {}{full} $\VB\sigma$-norm on the \lhs, 
we then obtain that
\Beq
  \frac \beta 2 \, \norma{\phi(t)}^2
  + \iot \norma{\phi(s)}_\Bs^2 \ds
  \,\leq\, \pier{c} \iot \norma{\phi(s)}^2 \ds
   \pier{{}+ \iot \bigl( \mu(s) , \phi(s) \bigr) \ds} .
  \label{diff2}
\Eeq
Finally, we add \eqref{dadiff13} to \eqref{diff2},
choose $\delta>0$ sufficiently small, and apply Gronwall's lemma.
We then conclude that $(\mu,\phi,S)=(0,0,0)$, and the proof is complete.

\Brem
In connection with Remark~\ref{Controls},
we could consider the equations obtained by adding the forcing terms, 
say, controls, to the \rhs s of the equations.
It is clear that no change is necessary in the above proof in order to obtain uniqueness 
also in this more general situation.
Furthermore, just minor modifications {}{lead to a continuous dependence result. More
precisely}, if $u_{\mu,i}$, $u_{\phi,i}$ and $u_{S,i}$, $i=1,2$, are two choices of the controls
and $u_\mu$, $u_\phi$ and $u_S$ denote their differences,
then we obtain, with the notation used in the proof,
\Bsist
  && \norma\mu_{\pier{\L2 H}}
  + \norma{1*\mu}_{\pier{\L\infty {\VA\rho}}}
  + \norma\phi_{\L\infty H\cap\L2{\VB\sigma}}
  + \norma S_{\L\infty H\cap\L2{\VC\tau}}
  \qquad
  \non
  \\
  && \leq\, \gianni{C_4}\, \bigl(
    \norma{u_\mu}_{\L2H}
    + \norma{u_\phi}_{\L2H}
    + \norma{u_S}_{\L2H}
  \bigr)\,,
  \label{dipcont}
\Esist
where \gianni{$C_4>0$} depends {}{only} on the structure, 
i.e., the linear operators, the shape of the nonlinearities,
the parameters $\alpha$ and~$\beta$, and the final time~$T$.
\Erem


\section{Existence}
\label{EXISTENCE}
\setcounter{equation}{0}

In this section, we prove the existence of a solution to problem \Pbl\ 
as stated in Theorem~\ref{Wellposedness}.
To help the reader, we start with a formal estimate {}{that
gives a flavor of the regularity to be expected
and, at the same time, indicates} the direction one can take for a rigorous proof.
Then, in the next subsections, we introduce the approximating problem and its discretization,
solve the discrete problem, perform rigorous estimates, 
and solve first the regularized problem and then problem~\Pbl.


\subsection{Preliminaries}
\label{PRELIMINARIES}
\setcounter{equation}{0}

Here is the formal estimate just mentioned.
We multiply \eqref{Iprima}, \eqref{Iseconda}, and \eqref{Iterza}, 
by~$\mu$, $-\dt\phi$, and~$S$, respectively, in the scalar product of~$H$.
Then we sum up and integrate over~$(0,t)$, where $t\in(0,T)$ is arbitrary,
noting that the terms involving the product $\mu\,\dt\phi$ cancel each other.
By accounting for \eqref{propA} and its analogues for the other two types of operators,
we obtain the identity
\Bsist
  && \frac \alpha 2 \, \norma{\mu(t)}^2 
  + \iot \norma{A^\rho\mu(s)}^2 \ds
  \non
  \\
  && \quad {}
  + \beta \intQt |\dt\phi|^2
  + \frac 12 \, \norma{B^\sigma\phi(t)}^2 
  + \iO F(\phi(t))
  \non
  \\
  \separa
  && \quad {}
  + \frac 12 \, \norma{S(t)}^2
  + \iot \norma{C^\tau S(s)}^2 \ds
  + \intQt P(\phi) (S-\mu)^2
  \non
  \\
  \separa
  && {} =\, \frac \alpha 2 \, \norma\muz^2 
  + \frac 12 \, \norma{B^\sigma\phiz}^2 
  + \iO F(\phiz)
  + \frac 12 \, \norma\Sz^2 \,.
  \non
\Esist
By reading $F_1$ instead of $\Ful$ in~\eqref{hpbelow}, 
and adding $\,|\Omega|C_0\,$ to both sides of the above equality,
we conclude that
\Bsist
  && \alpha^{1/2}\, \norma\mu_{\L\infty H}
  + \norma{A^\rho\mu}_{\L2H}
  \non
  \\[1mm]
  && \quad {}
  + \beta^{1/2} \,\norma{\dt\phi}_{\L2H}
  + \norma{B^\sigma\phi}_{\L\infty H}
  + \norma{F(\phi)+C_0}_{\L\infty\Luno}
  \non
  \\[1mm]
  && \quad {}
  + \norma S_{\L\infty H}
  + \norma{C^\tau S}_{\L2H}
  + \norma{P^{1/2}(\phi)(S-\mu)}_{\L2H}
  \non
  \\[1mm]
  && \leq\, C \,\bigl(
    \alpha^{1/2}\, \norma\muz
    + \norma{B^\sigma\phiz}
    + \norma{F(\phiz)+C_0}_1
    + \norma\Sz
  \bigr)\,,
  \label{formalestimate}
\Esist
where $C>0$ is a universal constant.
Eliminating~$C_0$ in the norms by means of the triangle inequality, we obtain an estimate 
\pier{that is nothing but \eqref{stimasoluz}.}

In the next subsections, after introducing and solving the discrete problem,
we implement the above \gianni{argument} to derive a rigorous a~priori estimate for the discrete solution.
Then, we use it for the necessary limiting procedures and solve the original problem.
For this purpose, it is convenient to introduce some notations at once.

\Bnot
\label{Interpolants}
Let $N$ be a positive integer, and \pier{let} $Z$ be one of the spaces $H$, $\VA\rho$, $\VB\sigma$,~$\VC\tau$.
We set $h:=T/N$ and $I_n:=((n-1)h,nh)$ for $n=1,\dots,N$.
Given $z=(z_0,z_1,\dots ,z_N)\in Z^{N+1}$,
we define the piecewise constant and piecewise linear interpolants
\Beq
  \overline z_h \in \L\infty Z , \quad
  \underline z_h \in \L\infty Z, 
  \aand
  \hat z_h \in \W{1,\infty}Z,
  \non
\Eeq
by setting 
\Bsist
  && \hskip -2em
  \overline z_h(t) = z^n
  \aand
  \underline z_h(t) = z^{n-1}
  \quad \hbox{for a.a.\ $t\in I_n$, \ $n=1,\dots,N$},
  \label{pwconstant}
  \\
  && \hskip -2em
  \hat z_h(0) = z_0
  \aand
  \dt\hat z_h(t) = \frac {z^{n+1}-z^n} h
  \quad \hbox{for a.a.\ $t\in I_n$, \ $n=1,\dots,N$}.
  \qquad
  \label{pwlinear}
\Esist
\Enot

For the reader's convenience,
we summarize the relations between the finite set of values
and the interpolants in the following proposition,
whose proof follows from \sfw\ computations.

\Bprop
\label{Propinterp}
With Notation~\ref{Interpolants}, we have that
\begin{align}
  & \norma{\overline z_h}_{\L\infty Z}
  = \max_{n=1,\dots,N} \norma{z^n}_Z \,, \quad
  \norma{\underline z_h}_{\L\infty Z}
   = \max_{n=0,\dots,N-1} \norma{z^n}_Z\,,
  \label{ouLinftyZ}
  \\
  & \norma{\dt\hat z_h}_{\L\infty Z}
  = \max_{0\leq n\leq N-1} \norma{(z^{n+1}-z^n)/h}_Z\,,
  \label{dtzLinftyZ}
  \\
  \separa
  & \norma{\overline z_h}_{\L2Z}^2
  = h \somma n1N \norma{z^n}_Z^2 \,, \quad
  \norma{\underline z_h}_{\L2Z}^2
  = h \somma n0{N-1} \norma{z^n}_Z^2 \,,
  \label{ouLdueZ}
  \\
  \separa
  & \norma{\dt\hat z_h}_{\L2Z}^2
  = h \somma n0{N-1} \norma{(z^{n+1}-z^n)/h}_Z^2\,, 
  \label{dtzLdueZ}
  \\
  & \norma{\hat z_h}_{\L\infty Z}
  = \max_{n=1,\dots,N}\, \max\,\{\norma{z^{n-1}}_Z,\norma{z^n}_Z\}
  = \max\,\{\norma{z_0}_Z,\norma{\overline z_h}_{\L\infty Z}\}\,,
  \label{hzLinftyZ}
  \\
  & \norma{\hat z_h}_{\L2Z}^2
  \leq h \somma n1N \bigl( \norma{z^{n-1}}_Z^2 + \norma{z^n}_Z^2 \bigr)
  \leq h\, \norma{z_0}_Z^2
  \,+\, 2\, \norma{\overline z_h}_{\L2Z}^2 \,.
  \label{hzLdueZ}
\end{align}
Moreover, it holds that
\begin{align}
  & \norma{\hat z_h(t)-\overline z_h(t)}_Z
  \leq \norma{\overline z_h(t)-\underline z_h(t)}_Z, \quad \
  \norma{\hat z_h(t)-\underline z_h(t)}_Z
  \leq \norma{\overline z_h(t)-\underline z_h(t)}_Z
  \non
  \\
  & \quad \aat \,,
  \label{diffae}
  \\[1mm]
  & \norma{\overline z_h-\hat z_h}_{\L\infty Z}
  = \max_{n=0,\dots,N-1} \norma{z^{n+1}-{z^n}}_Z
  = h \, \norma{\dt\hat z_h}_{\L\infty Z}\,,
  \label{diffLinfty}
  \\
  & \pier{\norma{\overline z_h-\hat z_h}^2_{\L\infty Z}
  \leq h \somma n0{N-1} h\left\|\frac{z^{n+1}-z^n}h\right\|_Z^2  
  = h \, \norma{\dt\hat z_h}_{\L2 Z}^2\,,}
  \label{pier2}
  \\
  & \norma{\overline z_h-\hat z_h}_{\L2Z}^2
  = \frac h3 \somma n0{N-1} \norma{z^{n+1}-z^n}_Z^2
  = \frac {h^2} 3 \, \norma{\dt\hat z_h}_{\L2Z}^2\,,
  \label{diffLdue}
\end{align}
and similar identities for the difference $\underline z_h-\hat z_h$.
As a consequence, we have the inequalities
\Bsist
  && \norma{\overline z_h-\underline z_h}_{\L\infty Z}
  \leq 2h \, \norma{\dt\hat z_h}_{\L\infty Z}\,,
  \qquad
  \label{diffbisLinfty}
  \\
  && \norma{\overline z_h-\underline z_h}_{\L2Z}^2
  \leq \frac {2 h^2} 3 \, \norma{\dt\hat z_h}_{\L2Z}^2 \,.
  \label{diffbisLdue}
\Esist
Finally, we have that
\Bsist
  && h \somma n0{N-1} \norma{(z^{n+1}-z^n)/h}_Z^2 
  \,\leq \,\norma{\dt z}_{\L2Z}^2\,,
  \non
  \\
  && \quad \hbox{if $z\in\H1Z$\aand $z^n=z(nh)$ for $n=0,\dots,N$}.
  \label{interpH1Z}
\Esist
\Eprop


\subsection{Approximation and discretization}
\label{APPROX}

In this subsection, we introduce an approximation of problem \Pbl\ 
and its time discretization.
Then, we solve the discrete problem.
We first introduce the Moreau--Yosida regularizations 
$\Ful$ and $\ful$ of $F_1$ of $f_1$ at the level $\lambda>0$
(see, e.g., \cite[p.~28 and p.~39]{Brezis}).
We set, for convenience,
\Beq
  \Fl := \Ful + F_2
  \aand
  \fl := \ful + f_2 \,.
  \label{yosida}
\Eeq
By~accounting for well-known properties of this regularization and the 
assumptions~\accorpa{hpFuno}{hpbelow}, we~have
\Beq
  \Ful(s) = \int_0^s \ful(s') \, ds' \,, \quad
  0 \leq \Ful(s) \leq F_1(s)\,,
  \aand
  \Fl(s) \geq - C_0\,,
  \label{propyosida}
\Eeq
for every $s\in\erre$,
as well~as
\Beq
  |\ful(s)| \leq |f_1^\circ(s)|
  \quad \hbox{for every $s\in D(f_1)$}\,,
  \label{propyosidabis}
\Eeq
where $f_1^\circ(s)$ is the element of $f_1(s)$ having minimum modulus.
By replacing $F_1$ in \eqref{seconda} by~$\Ful$, we obtain the following system:
\begin{align}
  & \alpha\, \< \dt\mul(t) , v >_\Ar
  + \bigl( \dt\phil(t) , v \bigr)
  + ( A^\rho \mul(t) , A^\rho v )
  = \bigl( \gianni{P(\phil(t)) (\Sl(t)-\mul(t))} , v \bigr)
  \non
  \\
  &\quad \hbox{for every $v\in\VA\rho$ and \aat},
  \label{primal}
  \\[1mm]
  \separa
  & \beta \,\bigl( \dt\phil(t) , \phil(t) - v \bigr) 
  + \bigl( B^\sigma\phil(t) , B^\sigma(\phil(t)-v) \bigr)
  \non
  \\
  & \quad {}
  + \iO \Ful(\phil(t))
  + \bigl( f_2(\phil(t)) , \phil(t)-v \bigr)
  \non
  \\
  & \leq \bigl( \mul(t) , \phil(t)-v \bigr)
  + \iO \Ful(v)
  \quad \hbox{for every $v\in\VB\sigma$ and \aat},
   \label{secondal}
  \\[1mm]
  \separa
  & \< \dt\Sl(t) , v >_\Ct
  + ( C^\tau\Sl(t) ,C^\tau v )
  = - \bigl( \gianni{P(\phil(t)) (\Sl(t)-\mul(t))} , v \bigr)
  \non
  \\
  & \quad \hbox{for every $v\in\VC\tau$ and \aat},
  \label{terzal}
  \\[1mm]
  & \mul(0) = \muz \,, \quad
  \phil(0) = \phiz\,, 
  \aand
  \Sl(0) = \Sz \,.
  \label{cauchyl}
\end{align}
\Accorpa\Pbll primal cauchyl
We stress that \eqref{secondal} is equivalent to both 
the time-integrated variational inequality
\Bsist
  && \beta \ioT \bigl( \dt\phil(t) , \phil(t) - v(t) \bigr) \, dt 
  + \ioT \bigl( B^\sigma\phil(t) , B^\sigma(\phil(t)-v(t)) \bigr) \, dt
  \non
  \\
  && \quad {}
  + \intQ \Ful(\phil)
  + \ioT \bigl( f_2(\phil(t)) , \phil(t)-v(t) \bigr) \, dt
  \non
  \\
  && \leq \ioT \bigl( \mul(t) , \phil(t)-v(t) \bigr) \, dt
  + \intQ \Ful(v)
  \quad \hbox{for every $v\in\L2{\VB\sigma}$},
  \label{intsecondal}
\Esist
and the pointwise variational equation 
(since $\Ful$ is differentiable and $\ful$ is its derivative)
\Bsist
  && \beta \bigl( \dt\phil(t) , v \bigr)
  + \bigl( B^\sigma\phil(t) , B^\sigma v \bigr)
  + \bigl( \fl(\phil(t)) ,  v \bigr)
  = \bigl( \mul(t) ,  v \bigr)
  \non
  \\[1mm]
  && \quad \hbox{for every $v\in\VB\sigma$ and \aat}.
  \label{eqsecondal}
\Esist

\Bthm
\label{Wellposednessl}
Under the same assumptions {}{as in} Theorem~\ref{Wellposedness}, the
problem \Pbll\ has \gianni{at least a} solution satisfying the analogues of \accorpa{regmu}{regS}.
\Ethm

\gianni{%
\Brem
\label{Uniqeps}
The above statement does not ensure uniqueness.
On the other hand, no uniqueness for the solution to the approximating problem is necessary for our purpose.
However, uniqueness is guaranteed if the spaces $\VA\rho$ and $\VC\tau$ satisfy~\eqref{embeddings}.
Indeed, in this case, what we have proved in Section~\ref{UNIQUENESS} can be applied
since $\Fl$ satisfies all the properties we have postulated for~$F$.
\Erem
}%

\gianni{The major part of the present section is devoted to the proof of Theorem~\ref{Wellposednessl}},
which is based on \pier{the discretization} procedure.
Thus, we introduce and solve the discrete problem and 
then take the limits of the interpolants as the time {}{step size} tends to zero.

\step
The discrete problem

We fix an integer $N>1$ and set $h:=T/N$.
Moreover, we fix a constant $L$ satisfying
\Beq
  L > \mathop{\rm Lip} f_2\,,
  \label{defL}
\Eeq
where $\mathop{\rm Lip} f_2$ is the \Lip\ constant of~$f_2$.
Then, the discrete problem consists in finding three $(N+1)$-tuples 
$(\mu^0,\dots,\mu^N)$, $(\phi^0,\dots,\phi^N)$, and $(S^0,\dots,S^N)$,
satisfying
\Bsist
  && \mu^0 = \muz \,, \quad
  \phi^0 = \phiz \,,
  \aand
  S^0 = \Sz   \,,
  \label{cauchyd}
  \\[2mm]
  && (\muu,\dots,\mu^N) \in (\VA{2\rho})^N\, , \quad
  (\phiu,\dots,\phi^N) \in (\VB{2\sigma})^N\,,
  \non
  \\
  && \aand
  (\Su,\dots,S^N) \in (\VC{2\tau})^N\,,
  \label{regd}
\Esist
and solving
\Bsist
  && \alpha \, \dhmun
  + \dhphin
  + A^{2\rho} \munp
  + P(\phin) \munp
  = P(\phin) \Snp\,,
  \qquad
  \label{primad}
  \\[1mm]
  && \beta \, \dhphin
  + B^{2\sigma} \phinp
  + (\fl + {}{L\,I})(\phinp)
  = L \phin + \munp \,,
  \label{secondad}
  \\
  && \dhSn 
  + C^{2\tau} \Snp
  + P(\phin) \Snp
  = P(\phin) \munp\,
  \label{terzad}
\Esist
\Accorpa\Pbld cauchyd terzad
\pier{a.e.~in~$\Omega$,} for $n=0,1,\dots,N-1$.
This problem can be solved inductively for $n=0,\dots,N-1$.
Namely, for a given $(\mun,\phin,\Sn)\in H\times H\times H$,
we show that there exists a unique triplet 
$(\munp,\phinp,\Snp)\in\VA{2\rho}\times\VB{2\sigma}\times\VC{2\tau}$
satisfying a problem equivalent to \accorpa{primad}{terzad}.
Here is the construction of the new problem.
We first observe that the linear operator 
\Beq
  \pier{ A_h v := {}{\frac \alpha h\,v} + A^{2\rho}v + P(\phin)v , \quad  v \in \VA{2\rho},}\label{defAh}
\Eeq
is an isomorphism \pier{from $\VA{2\rho}$ to $H$.} 
To see this, it suffices to apply \eqref{maxmon} to $A^{2\rho}$ and to notice that
the linear operator given by the last contribution $v\mapsto P(\phin)v$ 
is monotone and continuous from $H$ into itself, since $P$ is bounded and nonnegative.
By~the way, one also sees that $A_h^{-1}\in\calL(H;H)$ is monotone and that its norm is bounded by~$h/\alpha$.
Hence, \eqref{primad} can be solved for~$\munp$, and we can write
\Beq
  \munp 
  = A_h^{-1} \Bigl(
    \frac \alpha h \,\mun
    - \dhphin
    + P(\phin) \Snp
  \Bigr) .
  \label{primadbis}
\Eeq
So, we replace \eqref{primad} by \eqref{primadbis},
and \eqref{secondad} by the equation obtained by {}{inserting in \eqref{secondad} 
the expression \juerg{for} $\munp$ given by~\eqref{primadbis} in place of~$\munp$.
Thus, the new second equation reads}
\Bsist
  &&  \beta \, \dhphin
  + B^{2\sigma} \phinp
  + (\fl + {}{L\,I})(\phinp)
  \non
  \\
  && = \,L \phin 
  + A_h^{-1} \Bigl(
    \frac \alpha h \mun
    - \dhphin
    + P(\phin) \Snp
  \Bigr)
  \non
\Esist
or, even better,
\Bsist
  &&  \beta \, \dhphin
  + B^{2\sigma} \phinp
  + (\fl + {}{L\,I})(\phinp)
  + \frac 1h \, A_h^{-1} \phinp
  \non
  \\
  && = \,L \phin 
  + A_h^{-1} \Bigl(
    \frac \alpha h \mun
    + \frac 1h \, \phin
    + P(\phin) \Snp
  \Bigr)\,.
  \label{secondadbis}
\Esist
We rewrite \eqref{terzad} here, for convenience,
\Beq
  \dhSn 
  + C^{2\tau} \Snp
  + P(\phin) \Snp
  = P(\phin) \munp\,,
  \label{terzadbis}
\Eeq
and the new problem is given by \accorpa{primadbis}{terzadbis}.
We {}{now show that it can be solved} by a fixed point argument, provided that $h>0$ is small enough.
To this end, we construct some mappings.
In doing this, for simplicity, we use the symbols $\munp$, $\phinp$, and $\Snp$, 
as {}{if} they were independent variables.
The subscripts we choose remind the order of appearance of the equations.
Here are the mappings:
\Bsist
  & \Phi_3 : H \to \VC{2\tau} \subset H;
  & \quad \hbox{$\munp \mapsto \Snp$\quad by solving \eqref{terzadbis} for $\Snp$},
  \non
  \\
  & \Phi_2 : H \to \VB{2\sigma} \subset H;
  & \quad \hbox{$\Snp \mapsto \phinp$\quad by solving \eqref{secondadbis} for $\phinp$},
  \non
  \\
  & \Phi_1 : H \times H \to \VA{2\rho} \subset H;
  & \quad \hbox{$(\phinp,\Snp) \mapsto \munp$\quad by just applying \eqref{primadbis}},
  \non
  \\
  & \Phi : H \to H;
  & \quad 
  \overline\mu \mapsto
  \Phi_1 \Bigl(
    \Phi_2 \bigl(
      \Phi_3(\overline\mu)
    \bigr) 
    , \Phi_3(\overline\mu) 
  \Bigr) .
  \non
\Esist
Once we prove that these mappings are well defined and that $\Phi$ has a unique fixed point~$\mu^*$,
it is clear that the unique solution $(\munp,\phinp,\Snp)$ we are looking for is given by
$(\mu^*,\Phi_2(\Phi_3(\mu^*)),\Phi_3(\mu^*))$.
Let us start.
As for~$\Phi_3$, one adopts the same argument used to define~$A_h^{-1}$.
Concerning~$\Phi_2$, the proof is similar,
if one notes that the monotonicity of $A_h^{-1}$ follows from the one of~$A_h$
and that even $\fl$ and $f_2+{}{LI}$ are everywhere defined monotone operators, the last due to~\eqref{defL}.
Thus, all {}{of the mappings} are well defined.
Now, we consider~$\Phi_3$ and take any $\overline\mu_1\,,\,\overline\mu_2\in H$.
By writing \eqref{terzadbis} with $\overline S_i$ and~$\overline\mu_i$, $i=1,2$,
in place of $\Snp$ and~$\munp$, respectively,
and multiplying the difference by $\pier{\overline S_1 - \overline S_2}$, we immediately find~that
\Beq
  \frac 1h \, \norma{\overline S_1-\overline S_2}
  \leq \sup_{s\in\erre} P(s) \, \norma{\overline\mu_1-\overline\mu_2} \,.
  \non
\Eeq
This implies that
\Beq
  \norma{\Phi_3(\overline\mu_1) - \Phi_3(\overline\mu_2)}
  \leq K_3 h \, \norma{\overline\mu_1 - \overline\mu_2}
  \quad \hbox{for every $\overline\mu_1,\overline\mu_2\in H$}\,,
  \non
\Eeq 
where $K_3$ is the supremum of~$P$.
Similarly, one {}{shows} that
\Bsist
  & \norma{\Phi_2(\overline S_1)-\Phi_2(\overline S_2)}
  \leq K_2 \, h \, \norma{\overline S_1-\overline S_2}\,,
  \non
  \\[1.5mm]
  & \norma{\Phi_1(\overline\phi_1,\overline S_1) - \Phi_1(\overline\phi_2,\overline S_2)}
  \leq K_1 \bigl(
  \norma{\overline\phi_1 - \overline\phi_2}
  + \norma{\overline S_1 - \overline S_2}
  \bigr)\,,
  \non
\Esist
for every $\overline S_i\in H$ and $\overline\phi_i\in H$, $i=1,2$,
and some constants $K_2$ and~$K_1$.
Hence, \pier{there is a constant~$K$ such that}
\Beq
  \norma{\Phi(\overline\mu_1) - \Phi(\overline\mu_2)}
  \leq K h \, \norma{\overline\mu_1 - \overline\mu_2}
  \quad \hbox{for every $\overline\mu_1,\overline\mu_2\in H$}.
  \non
\Eeq
Therefore, if $Kh<1$, $\Phi$ is a contraction in $H$ and thus has a unique fixed point.
We conclude that the discrete problem is uniquely solvable by assuming that {}{$0<h<K^{-1}$.}


\subsection{Solution of the approximating problem}
\label{SOLUZAPPROX}

As announced in the Introduction,
we prove the existence of a solution to the approximating problem~\Pbll\
by taking the limit of the interpolants of the solution to the discrete problem 
as the time step size $h$ tends to zero.
According to Notation~\ref{Interpolants},
we remark at once that the regularity required for discrete solution implies that 
\Bsist
  && \muh \in \L\infty{\VA\rho} , \quad
  \undermuh \in \pier{\L\infty{H}},
  \aand
  \overmuh \in \L\infty{\VA{2\rho}},
  \label{regmuh}
  \\[1mm]
  && \phih \in \W{1,\infty}{\VB\sigma} , \quad
  \underphih \in \L\infty{\VB\sigma},
  \aand
  \overphih \in \L\infty{\VB{2\sigma}},
  \qquad
  \label{regphih}
  \\[1mm]
  && \Sh \in \L\infty{\VC\tau} , \quad
  \underSh \in \pier{\L\infty{H}},
  \aand
  \overSh \in \L\infty{\VC{2\tau}},
  \label{regSh}
\Esist
and that the discrete problem also reads
\Bsist
  && \alpha \, \dt\muh
  + \dt\phih
  + A^{2\rho} \overmuh
  + P(\underphih) \overmuh
  = P(\underphih) \overSh
  \quad \pier{\aeQ},
  \label{primah}
  \\[1mm]
  && \beta \, \dt\phih
  + B^{2\sigma} \overphih
  + ( \fl + {}{L\,I} ) (\overphih)
  = L \underphih
  + \overmuh
  \quad \pier{\aeQ} ,
  \label{secondah}
  \\[1mm]
  && \dt\Sh
  + C^{2\tau} \overSh
  + P(\underphih) \overSh
  = P(\underphih) \overmuh
  \quad \pier{\aeQ} ,
  \label{terzah}
  \\[1mm]
  && \muh(0) = \muz \,, \quad
  \phih(0) = \phiz\,, \quad
  \Sh(0) = \Sz \, \quad  \pier{\aeO} .
  \label{cauchyh}
\Esist
\Accorpa\Pblh primah cauchyh
\pier{We  point out that \juergen{the} equations~\eqref{primah}--\eqref{terzah} have been written 
a.e.\ in $Q$, and in this case all \juergen{of the} terms, including 
$A^{2\rho} \overmuh$, $B^{2\sigma} \overphih$, and  $C^{2\tau} \overSh$,
are interpreted as functions of space and time; another way of reading~\eqref{primah}--\eqref{terzah} could be in $H$, a.e.\ in $(0,T)$, as the single terms make sense in the space $H$ as well.}

So, our aim is \juerg{to let} $h$ tend to zero in \pier{\Pblh}
(or~in some equivalent formulation).
Hence, we start estimating.
We do this on the solution to the discrete problem~\Pbld,
by adapting the procedure that led to the formal estimate of Section~\ref{PRELIMINARIES}.
Then, we express the bounds we find in terms of the interpolants.
According to the general rule stated at the end of Section~\ref{STATEMENT},
the (possibly different) values of the constants termed $c$ 
are independent of the parameters $h$, $\lambda$, $\alpha$, and~$\beta$.

\step
Basic a priori estimate

We test \eqref{primad}, \eqref{secondad} and~\eqref{terzad}
(by~taking the scalar product in~$H$)
by~$\munp$, $(\phinp-\phin)/h$ and $\Snp$, respectively,
and add the resulting identities to each other.
Noting an obvious cancellation, we obtain the equality
\Bsist
  && \frac\alpha h \, (\munp , \munp-\mun)
  + \norma{A^\rho \munp}^2
  + \iO P(\phin) (\munp-\Snp)^2
  \non
  \\
  && \quad {}
  + \beta \, \Norma{\dhphin}^2
  + \frac 1h \, \bigl( B^\sigma\phinp , B^\sigma(\phinp-\phin) \bigr)
  \non
  \\
  \separa
  && \quad {}
  + \frac 1h \, \bigl( (\fl + L)(\phinp) , \phinp-\phin \bigr)
  \non
  \\
  && \quad {}
  + \frac 1h \, (\Snp , \Snp-\Sn)
  + \norma{C^\tau \Snp}^2
  \non
  \\
  \separa
  && = \, \frac Lh \, ( \phin , \phinp-\phin ) .
  \non
\Esist
Now, we observe that the function $s\mapsto\Fl(s)+\frac L2\,r^2=\Ful(s)+F_2(s)+\frac L2\,r^2$
is convex on~$\erre$, since $\Ful$ is convex and $L$ satisfies~\eqref{defL}.
Thus, we have that
\Bsist
  && \bigl( (\fl + {}{L\,I})(\phinp) , \phinp-\phin \bigr)
  \non
  \\  
  && \geq \iO \Fl(\phinp) + \frac L2 \, \norma\phinp^2 
  - \iO \Fl(\phin) - \frac L2 \, \norma\phin^2 .
  \non
\Esist
Therefore, by using this inequality and applying the identity \eqref{elementare} 
to some terms of the previous equality, we deduce that
\Bsist
  && \frac \alpha {2h} \, \norma\munp^2
  + \frac \alpha {2h} \, \norma{\munp-\mun}^2
  - \frac \alpha {2h} \, \norma\mun^2
  + \norma{A^\rho\munp}^2
  \non
  \\
  && \quad {}
  + \iO P(\phin) (\munp-\Snp)^2
  + \beta \, \Norma{\dhphin}^2
  \non
  \\
  && \quad {}
  + \frac 1{2h} \, \norma{B^\sigma \phinp}^2
  + \frac 1{2h} \, \norma{B^\sigma(\phinp-\phin)}^2
  - \frac 1{2h} \, \norma{B^\sigma\phin}^2
  \non
  \\
  \separa
  && \quad {}
  + \frac 1h \iO \Fl(\phinp) + \frac L{2h} \, \norma\phinp^2 
  - \frac 1h \iO \Fl(\phin) - \frac L{2h} \, \norma\phin^2 
  \non
  \\
  \separa
  && \quad {}
  + \frac 1{2h} \, \norma\Snp^2
  + \frac 1{2h} \, \norma{\Snp-\Sn}^2
  - \frac 1{2h} \, \norma\Sn^2
  + \norma{C^\tau\Snp}^2
  \non
  \\
  \separa
  && \leq \,\frac L{2h} \, \norma\phinp^2
  - \frac L{2h} \, \norma\phin^2
  - \frac L{2h} \, \norma{\phinp-\phin}^2 .
  \non
\Esist
At this point, we first note two cancellations; 
then, we multiply by $h$ and sum up with respect to $n=0,\dots,m-1$
for $m=1,\dots,N$.
We obtain
\Bsist
  && \frac \alpha 2 \, \norma\mum^2
  - \frac \alpha 2 \, \norma\muz^2
  + \frac \alpha 2 \, \somma n0{m-1} \norma{\munp-\mun}^2
  + \somma n0{m-1} h\, \norma{A^\rho\munp}^2
  \non
  \\
  && \quad {}
  + \somma n0{m-1} h \iO P(\phin) (\munp-\Snp)^2
  + \beta \somma n0{m-1} h \, \Norma{\dhphin}^2
  \non
  \\
  && \quad {}
  + \frac 12 \, \norma{B^\sigma \phim}^2
  - \frac 12 \, \norma{B^\sigma\phiz}^2
  + \frac 12 \, \somma n0{m-1} \norma{B^\sigma(\phinp-\phin)}^2
  + \iO \Fl(\phim) 
  - \iO \Fl(\phiz) 
  \non
  \\
  \separa
  && \quad {}
  + \frac 12 \, \norma\Sm^2
  - \frac 12 \, \norma\Sz^2
  + \frac 12 \, \somma n0{m-1} \norma{\Snp-\Sn}^2
  + \somma n0{m-1} h \, \norma{C^\tau\Snp}^2
  \non
  \\
  \separa
  && \leq \,-\, \frac L2 \, \somma n0{m-1} \norma{\phinp-\phin}^2 .
  \non
\Esist
Clearly, this inequality also holds for $m=0$
if it is understood that all the sums vanish since the set of the indices is empty.
Therefore, by rearranging, 
accounting for \eqref{propyosida}, adding $|\Omega|C_0$ to both sides 
and owing to the assumption \eqref{hpdati} on the initial data,
we obtain an estimate (the~analogue of~\eqref{formalestimate})
that in terms of the interpolants reads
\Bsist
  && \alpha \,\norma\overmuh_{\L\infty H}^2
  + \frac \alpha h \, \norma{\overmuh-\undermuh}_{\L2H}^2
  + \norma{A^\rho \overmuh}_{\L2H}^2
  \non
  \\
  && \quad {}
  + \norma{(P(\underphih))^{1/2}(\overmuh-\overSh)}_{\pier{\L2H}}^2  
  + \beta \, \norma{\dt\phih}_{\L2H}^2
  + \pier{\frac Lh } \, \norma{\overphih-\underphih}_{\L2H}^2
  \non
  \\
  && \quad {}
  + \norma{B^\sigma\overphih}_{\L\infty H}^2
  + \frac 1h \, \norma{B^\sigma(\overphih-\underphih)}_{\L2H}^2
  + \norma{\Fl(\overphih)}_{\L\infty\Luno}
  \non
  \\
  \separa
  && \quad {}
  + \norma\overSh_{\L\infty H}^2
  + \frac 1h \, \norma{\overSh-\underSh}_{\L2H}^2
  + \norma{C^\tau \overSh}_{\L2H} ^2
  \non
  \\
  && \leq \,C_0'  \,\bigl(
    \alpha \, \norma\muz^2
    + \norma{B^\sigma\phiz}^2
    + \norma{F(\phiz)}_1
    + \norma\Sz^2
    + 1
  \bigr) \,,
  \label{basic}
\Esist
where $\,C_0'$\, depends only on $\Omega$ and the constant~$C_0$.

\step
First consequences

We observe that (see also \pier{\eqref{pier2}})
\begin{align*}
  &\norma{\overphih(t)}
  \leq \norma{\phih(t)}
  + \norma{\overphih(t)-\phih(t)}
 \\
 &\leq \norma\phiz + T^{1/2} \norma{\dt\phih}_{\L2H}
  + \pier{h^{1/2} \norma{\dt\phih}_{\L2H}}
  \non
\end{align*}
\aat.	
Moreover, the inequality $P(s)\leq c\,(P(s))^{1/2}$ holds true for every $s\in\erre$
due to the boundedness of~$P$.
Hence, we infer from \eqref{basic} that
\Bsist
  && \norma{\overmuh}_{\L\infty H\cap\L2{\VA\rho}}
  + \norma{\overphih}_{\pier{\L\infty{\VB\sigma}}}
  + \norma{\overSh}_{\L\infty H\cap\L2{\VC\tau}}
  \non
  \\
  && \quad {}
  + \norma\phih_{\H1H}
  + \norma{\Fl(\overphih)}_{\L\infty\Luno}
  + \norma{P(\underphih)(\overmuh-\overSh)}_{\L2H}
  \non
  \\
  && \leq c_{\alpha,\beta}\,,
  \label{primastimah}
\Esist
as well as (\pier{due to \eqref{diffae}})
\Bsist
  && \norma{\overmuh-\muh}_{\L2H}
  + \norma{\overphih-\underphih}_{\L2{\VB\sigma}}
  + \norma{\overphih-\phih}_{\L2H}
  \non
  \\
  && \quad {}
  + \norma{\overSh-\Sh}_{\L2H}
  \leq c_\alpha \, h^{1/2} .
  \label{stimediffh}
\Esist
By combining with \eqref{primastimah}, we deduce that
\Beq
  \norma\underphih_{\L2{\VB\sigma}} \leq c_{\alpha,\beta} \,.
  \label{stimaunderphih}
\Eeq
We also derive an estimate that we will use later on.
Since $F_2$ grows at most quadratically due to~\eqref{hpFdue},
the inequality \eqref{primastimah} yields an estimate for $F_2(\overphih)$ 
in~$\pier{\L\infty{L^1(\Omega)}}$.
Therefore, owing to the estimate of $\Fl(\overphih)$ given by~\eqref{basic},
we deduce~that
\Beq
  \norma{\Ful(\overphih)}_{\pier{\L\infty{L^1(\Omega)}}} \leq c_{\alpha,\beta} \,.
  \label{stimaFh}
\Eeq

\step
Second a priori estimate

By direct computation, for $n=0,\dots,N-1$ and for a.e.\ $t\in(nh,(n+1)h)$, we have that
\Bsist
  && \norma{B^\sigma(\phih(t)-\underphih(t))}
  = \norma{B^\sigma(\phin+\textstyle\frac{t-nh}h\,(\phinp-\phin)-\phin)}
  \non
  \\
  && = \textstyle\frac{t-nh}h \, \norma{B^\sigma(\phinp-\phin)}
  = \textstyle\frac{t-nh}h \, \norma{B^\sigma(\overphih(t)-\underphih(t))}
  \leq \norma{B^\sigma(\overphih(t)-\underphih(t))}\,,
  \non
\Esist
whence
\Beq
  \norma{B^\sigma(\phih-\underphih)}_{\L2H}
  \leq \norma{B^\sigma(\overphih-\underphih)}_{\L2H} \,.
  \non
\Eeq
By also accounting for \eqref{stimediffh}, we deduce that
\Beq
  \norma{\phih-\underphih}_{\L2{\VB\sigma}}
  \leq c_\alpha \, h^{1/2} \,,
  \non
\Eeq
and \eqref{stimaunderphih} yields that
\Beq
  \norma\phih_{\L2{\VB\sigma}}
  \leq c_{\alpha,\beta} \,.
  \label{secondastimah}
\Eeq

\step
Third a priori estimate

By equation \eqref{primah} and assumption \eqref{hpP}, we have
\Bsist
  && \alpha \, \norma{\dt\muh}_{\L2{\VA{-\rho}}}
  \non
  \\[1mm]
  && \leq c \, \bigl(
    \norma{\dt\phih}_{\L2H}
    + \norma{A^{2\rho}\overmuh}_{\L2{\VA{-\rho}}}
    + \norma\overmuh_{\L2H}
    + \norma\overSh_{\L2H}
  \bigr) .
  \non 
\Esist
Then, we account for \eqref{primastimah} and the first of \eqref{extensions2}
to obtain an estimate for the time derivative~$\dt\muh$.
By proceeding analogously with equation~\eqref{terzah}, we conclude that
\Beq
  \norma{\dt\muh}_{\L2{\VA{-\rho}}}
  + \norma{\dt\Sh}_{\L2{\VC{-\tau}}}
  \leq c_{\alpha,\beta} \,.
  \label{terzastimah}
\Eeq

\step
Convergence 

By recalling \accorpa{primastimah}{terzastimah},
we see that there exist a triplet $(\mul,\phil,\Sl)$ such~that
\Bsist
  & \overmuh \to \mul
  & \quad \hbox{weakly star in $\L\infty H\cap\L2{\VA\rho}$}\,,
  \label{convovermuh}
  \\
  & \muh \to \mul
  & \quad \hbox{weakly star in \pier{$\H1{\VA{-\rho}}\cap\L\infty H$}}\,,
  \label{convmuh}
  \\
  \separa
  & \overphih \to \phil
  & \quad \hbox{weakly \pier{star in $\L\infty{\VB\sigma}$}}\,,
  \label{convoverphih}
  \\
  & \underphih \to \phil
  & \quad \hbox{weakly in $\L2{\VB\sigma}$}\,,
  \label{convunderphih}
  \\
  & \phih \to \phil
  & \quad \hbox{weakly in $\H1H\cap\L2{\VB\sigma}$}\,,
  \label{convphih}
  \\
  \separa
  & \overSh \to \Sl
  & \quad \hbox{weakly star in $\L\infty H\cap\L2{\VC\tau}$}\,,
  \label{convoverSh}
  \\
  & \Sh \to \Sl
  & \quad \hbox{weakly star in $\pier{\H1{\VC{-\tau}}\cap\L\infty H}$}\,,
  \label{convSh}
\Esist
at least for some sequence $h_k\searrow0$.
From \eqref{convmuh}, \eqref{convphih}, \eqref{convSh}, and~\eqref{cauchyh},
we deduce that the initial conditions \eqref{cauchyl} are satisfied by the limiting triplet.
Next, we prove that \accorpa{primal}{terzal} are fulfilled as well.
By first applying the Aubin--Lions lemma
(see, e.g., \cite[Thm.~5.1, p.~58]{Lions}) to $\phih$ on account of \pier{\eqref{convphih}}, and then owing to \eqref{stimediffh},
we deduce that
\Beq
  \phih \to \phil , \quad
  \overphih \to \phil,
  \aand
  \underphih \to \phil,
  \quad \hbox{strongly in $\L2H$}.
  \label{strongphih}
\Eeq
\pier{In particular note that the limit $\phil$ is in $\H1H\cap\L\infty{\VB\sigma}$,
thanks to \eqref{convoverphih} and \eqref{convphih}. 
Next, by} recalling that $f_2$ and $P$ are \Lip\ continuous
(see \eqref{hpFdue}, \eqref{hpP}, and~\eqref{deffunodue}),
and that the same holds for~$\ful$ due to the general properties of the Yosida approximation,
we infer that
\Beq
  \fl(\overphih) \to \fl(\phil)
  \aand
  P(\underphih) \to P(\phil)
  \quad \hbox{strongly in $\L2H$}.
  \non
\Eeq
The latter, \eqref{convoverSh}, and \eqref{convovermuh} imply that
\Beq
  P(\underphih) (\overSh - \overmuh)  \to P(\phil) (\Sl - \mul)
  \quad \hbox{weakly in $\LQ1$} . 
  \non
\Eeq
On the other hand, $P(\underphih) (\overSh - \overmuh)$ is bounded in $\L2H$ by~\eqref{primastimah}.
Therefore, we conclude that
\Beq
  P(\underphih) (\overSh - \overmuh)  \to P(\phil) (\Sl - \mul)
  \quad \hbox{weakly in $\L2H$} . 
  \non
\Eeq
\pier{In view of \eqref{strongphih}, we have that, possibly taking another subsequence of $h$, 
\Beq
  \overphih (t) \to \phil (t)
    \quad \hbox{strongly in $H$}, \quad \aat .
  \non
\Eeq 
Hence, by lower semicontinuity it turns out that 
\Beq
  \iO\Ful(\phil (t) )
  \leq \liminf_{h\searrow0} \iO \Ful(\overphih (t))
  \leq c_{\alpha,\beta} \quad \aat . 
  \label{stimaFl}
\Eeq}%
At this point, we write \accorpa{primah}{terzah} in the equivalent form
\Bsist
  && \ioT \Bigl(
    \alpha \, \< \dt\muh(s) , v(s) >_{A,\rho}
    + ( \dt\phih(s) , v(s))
    + ( A^\rho\overmuh(s) , A^\rho v(s) )
  \Bigr) \ds
  \non
  \\
  && = \ioT \bigl( P(\underphih(s))  (\overSh(s) - \overmuh(s)) , v(s) \bigr) \ds
  \qquad \hbox{for every $v\in\L2{\VA\rho}$}\,,
  \non
  \\[1mm]
  && \ioT \Bigl(
    \bigl( \beta \, \dt\phih(s) , v(s) \bigl)
    + \bigl( B^\sigma \overphih(s) , B^\sigma v(s) \bigr)
    + \bigl( ( \ful + f_2 + {}{L\,I} ) (\overphih(s)) , v(s) \bigr)
  \Bigr) \ds
  \non
  \\
  && = \ioT \bigl( L \underphih(s) + \overmuh(s), v(s) \bigr) \ds
  \qquad \hbox{for every $v\in\L2{\VB\sigma}$}\,,
  \non
  \\[1mm]
  && \ioT \Bigl(
    \bigl( \dt\Sh(s) , v(s) \bigr)
    + \bigl( C^\tau \overSh(s){}{,C^\tau} v(s) \bigr)
  \Bigr) \ds
  \non
  \\
  && = - \ioT \bigl( P(\underphih(s))  (\overSh(s) - \overmuh(s)) , v(s) \bigr) \ds
  \qquad \hbox{for every $v\in\L2{\VC\tau}$}\,,
  \non
\Esist
and let $h$ tend to zero on account of the convergence properties we have established.
We obtain the integrated versions of \eqref{primal}, \eqref{terzal}, and \eqref{eqsecondal}. \pier{Now, starting from \eqref{eqsecondal}, 
we can perform the formal procedure that led to the estimate~\eqref{formalestimate},
by observing that the argument used there is now correct.
One obtains the estimate}
\begin{align}
  & \alpha^{1/2}\, \norma\mul_{\L\infty H}
  + \norma{A^\rho\mul}_{\L2H}
  \non
  \\[1mm]
  & \quad {}
  + \beta^{1/2} \,\norma{\dt\phil}_{\L2H}
  + \norma{B^\sigma\phil}_{\L\infty H}
  + \norma{\Fl(\phil)+C_0}_{\L\infty\Luno}
  \non
  \\[1mm]
  & \quad {}
  + \norma\Sl_{\L\infty H}
  + \norma{C^\tau\Sl}_{\L2H}
  + \norma{P^{1/2}(\phil)(\Sl-\mul)}_{\L2H}
  \non
  \\[1mm]
  & \leq C \bigl(
    \alpha^{1/2} \norma\muz
    + \norma{B^\sigma\phiz}
    + \norma{\Fl(\phiz)+ C_0}_1
    + \norma\Sz
  \bigr),
  \label{truestimate}
\end{align}
\gianni{where $C_0$ is given by~\eqref{propyosida} and $C$ is a universal constant.}
\pier{Just something on the regularity is missing, namely, 
the requirements for the time derivatives $\dt\mul$ and~$\dt\Sl$. 
But these regularities immediately follow} from~\eqref{terzastimah}, which also yields \juergen{that}
\Beq
  \norma{\dt\mul}_{\L2{\VA{-\rho}}}
  + \norma{\dt\Sl}_{\L2{\VC{-\tau}}}
  \leq c_{\alpha,\beta} \,.
  \label{stamdtmulSl}
\Eeq
This concludes the proof of Theorem~\ref{Wellposednessl}.

\subsection{Solution to the original problem}
\label{SOLUZ}

In this section, we conclude the proof of Theorem~\ref{Wellposedness}.
Namely, we costruct a solution $\soluz$ by letting $\lambda$ tend to zero in the approximating problem.
\gianni{From \accorpa{truestimate}{stamdtmulSl} and the boundedness of~$P$
(which implies $P\leq c\,P^{1/2}$), we derive the following estimate:}
\Bsist
  && \norma\mul_{\H1{\VA{-\rho}}\cap\L\infty H\cap\L2{\VA\rho}}
  + \norma\phil_{\H1H\cap\gianni{\L\infty{\VB\sigma}}}
  \non
  \\[1mm]
  &&
  + \,\norma\Sl_{\H1{\VC{-\tau}}\cap\L\infty H\cap\L2{\VC\tau}}
  + \norma{P(\phil)(\mul-\Sl)}_{\L2H}
  \,\leq\, c_{\alpha,\beta} \,.
  \qquad
  \label{stimal}
\Esist
Therefore, by using the same arguments of the previous subsection \pier{and the generalized Ascoli theorem,} we deduce that (for some sequence $\lambda_k\searrow0$)
\Bsist
  & \mul \to \mu
  & \quad \hbox{weakly in $\H1{\VA{-\rho}}\cap\L2{\VA\rho}$ and strongly in $\L2H$},
  \non 
  \\
  & \phil \to \phi
  & \quad \hbox{weakly star in $\H1H\cap\L\infty{\VB\sigma}$ and strongly in \pier{$\C0H$}},
  \non 
  \\
  & \Sl \to S
  & \quad \hbox{weakly in $\H1{\VC{-\tau}}\cap\L2{\VC\tau}$ and strongly in $\L2H$} .
  \non 
\Esist
Similarly as before, we obtain the initial conditions, and we also have that
\Bsist
  && f_2(\phil) \to f_2(\phi)
  \aand
  P(\phil) \to P(\phi)
  \quad \hbox{strongly in $\L2H$},
  \non
  \\
  && P(\phil)(\Sl-\mul) \to P(\phi)(S-\mu)
  \quad \hbox{weakly in $\L2H$ \pier{and strongly in $L^1(Q)$}} .\qquad
  \non
\Esist
In particular, we can pass to the limit in \eqref{primal} and \eqref{terzal} to
 obtain \eqref{prima} and~\eqref{terza}, respectively.
On the contrary, some more work has to be done for the equation for~$\phi$,
\pier{in particular to argue on the $\liminf$ in the \lhs\ of the inequality~\eqref{secondal} or, equivalently, \eqref{intsecondal}. 
To this concern, we can show that
\Beq
\iO F_1(\phi (t) ) 
  \leq \liminf_{\lambda\searrow0} \iO \Ful (\phil (t) ) \quad \hbox{for all } t\in [0,T]. 
\label{pier3} 
\Eeq
Indeed, let us recall the definitions of the resolvent $J_\lambda$ 
of $f_1 = \partial F_1$ and the Moreau-Yosida approximation $\Ful$ of $F_1$, \juergen{which are given by}
$$
  J_\lambda:=(I+\lambda f_1)^{-1} , \quad \ 
  \Ful(r) := \min_{s\in\erre} \bigl\{ {\textstyle\frac 1{2\lambda}} \, |s-r|^2 + F_1(s) \bigr\},
$$
in order to point out the property (see, e.g., \cite[Prop.~2.11, p.~39]{Brezis}) 
\Beq
  \Ful(r) = F_1 (J_\lambda(r)) +\frac 1{2\lambda} |J_\lambda(r) - r|^2 
  \quad \hbox{for all $r\in\erre $}.
  \label{formula}
\Eeq
Now, we know that $\phil (t)$ converges to $ \phi (t)$ in $H$ as $\lambda \searrow 0$ and that 
$\iO \Ful (\phil (t) )$ is nonnegative and bounded independently of $\lambda$, by virtue of \eqref{truestimate} and \eqref{hpbelow}. Then, using the representation \eqref{formula}, it \gianni{immediately follows} that 
$\iO F_1 (J_\lambda \phil (t)) $ is bounded independently of $\lambda$ and 
that also $J_\lambda \phil (t)$ converges to $ \phi (t)$ in $H$ as $\lambda \searrow 0$.
Hence, \eqref{pier3} follows from the lower semicontinuity of the convex functional 
$v\mapsto \iO F_1 (v) $ in $H$. In addition, this argument also entails that 
$F_1 (\phi) \in \L\infty{L^1(\Omega)} $ and 
\Beq
0 \leq \intQ F_1(\phi) 
  \leq \liminf_{\lambda\searrow0} \intQ \Ful (\phil ) , 
\label{pier4} 
\Eeq
which ensures \eqref{regFuphi}. At this point,
since 
$$
\ioT \bigl( B^\sigma\phi(t) , B^\sigma \phi(t)) \bigr) \, dt
\leq \liminf_{\lambda\searrow0}
\ioT \bigl( B^\sigma\phil(t) , B^\sigma \phil(t)) \bigr) \, dt
$$
by \juergen{the} weak lower semicontinuity of the norm in $\L2H$, it is sufficient 
to let $\lambda$ tend to zero in \eqref{intsecondal} in order to obtain~\eqref{intseconda}.} \gianni{This concludes the proof of the existence of a solution in the sense of Theorem~\ref{Wellposedness};
it remains to complete the proof of the estimates 
\accorpa{stimasoluz}{dapier} for the solution we have constructed.}

\pier{In view of \eqref{truestimate} and \eqref{hpbelow}, we claim that
\begin{align}
  0 \leq \iO (F(\phi (t) ) + C_0) 
  &\leq \liminf_{\lambda\searrow0} \iO \left(\Ful (\phil (t) ) + F_2 (\phil (t)) + C_0\right)
  \non\\
  &\leq C \bigl(
    \alpha^{1/2} \norma\muz
    + \norma{B^\sigma\phiz}
    + \norma{F(\phiz)+ C_0}_1
    + \norma\Sz
  \bigr)
  \label{forF1}
\end{align}
for all $t\in [0,T]$. Indeed, the first inequality in \eqref{forF1} is a consequence of \eqref{hpbelow} when taking the limit as  $\lambda \searrow 0$. Moreover,  as, for all $t\in [0,T]$, $\phil (t)$ converges to $ \phi (t)$ in $H$ and $F_2 \in C^1(\erre) $ has a Lipschitz continuous derivative (i.e., $f_2$), using the Taylor formula it is not difficult to verify that
$$  F_2(\phil(t) ) \to F_2(\phi (t))
  \quad \hbox{strongly in $L^1(\Omega)$.} 
$$
Furthermore, the last term in \eqref{forF1} \gianni{comes from the \rhs\ of \eqref{truestimate}
as a consequence of}  $0 \leq \iO \Ful (\phiz) \leq \iO F_1(\phiz)$, \gianni{and $\iO F_1(\phiz)$} if finite because of \eqref{hpdati}. Then, \eqref{forF1} follows easily from \eqref{pier3}.}

\pier{%
Now, by \eqref{truestimate}, \eqref{forF1}, 
and the weak or weak star lower semicontinuity of norms, we easily obtain
 \eqref{stimasoluz}.}
\gianni{%
As for~\eqref{stimadtsoluz}, we find the right bound 
for the time derivative $\dt(\alpha\mu+\phi)$.
\juergen{We} observe that \eqref{primal} yields for every $v\in\L2{\VA\rho}$ \juergen{that}
\begin{align}
  & \ioT \< \dt( \alpha \mul + \phil )(t), v(t) >_\Ar \, dt
  \non
  \\
  & = - \ioT \bigl( A^\rho \mul(t) , A^\rho v(t) \bigr) \, dt
  + \ioT \bigl( P(\phil(t)) (\Sl(t) - \mul(t)) , v(t) \bigr) \, dt
  \non
  \\
  & \leq \norma{A^\rho\mul}_{\L2H} \, \norma{A^\rho v}_{\L2H}
  + \pier{\left(\sup P^{1/2} \right)} \norma{P^{1/2}(\phil)(\Sl-\mul)}_{\L2H} 
  \, \norma v_{\L2H}
  \non
  \\[1mm]
  & \leq \bigl\{
    \norma{A^\rho\mul}_{\L2H}
    + \pier{\left(\sup P^{1/2} \right)} \norma{P^{1/2}(\phil)(\Sl-\mul)}_{\L2H}
  \bigr\} \norma v_{\L2{\VA\rho}} \,.
  \non
\end{align}
This, \pier{along with \eqref{truestimate}, provides} the analogue of the desired estimate for $\dt(\alpha\mul+\phil)$,
and the estimate for $\dt(\alpha\mu+\phi)$ follows immediately.
Since the treatment of $\dt S$ is quite similar, \pier{\eqref{stimadtsoluz}} is completely proved.
Finally, to obtain \eqref{dapier},
it suffices to remark that the further assumption~\eqref{pier} we make implies~that
\pier{%
\Bsist
  && \norma\phi_{\L\infty{\VB\sigma}}^2
  = \norma\phi_{\L\infty H}^2
  + \norma{B^\sigma \phi}_{\L\infty H}^2
  \non
  \\
  && \leq \frac 1 {c_1} \, \bigl( \norma{F(\phi)}_{\L\infty\Luno} + c_2 \bigr)
  + \norma{B^\sigma \phi}_{\L\infty H}^2\,,
  \non
\Esist%
so that \eqref{stimasoluz} plainly leads to the correct estimate \eqref{dapier}.
Then, Theorem~\ref{Wellposedness} turns out to be completely proved.}}


\section{Regularity}
\label{REGULARITY}
\setcounter{equation}{0}

\gianni{This section is devoted to establish further properties of the solution to problem \Pbl.
Namely, we prove Theorems~\ref{Regularity} and \ref{Secondeq} as well as Corollary~\ref{Fullregularity}.
We start with the first of these results.}

\step
Proof of Theorem~\ref{Regularity}

\gianni{The rigorous proof is based on a priori estimates for the solution to the discrete problem
obtained by first performing the discrete differentiation of \eqref{secondad} 
and then suitably testing the resulting equality as well as \eqref{primad} and~{}{\eqref{terzad}}, 
and finally summing up.
Since the details are rather heavy, we prefer to deal with the approximating problem \Pbll, directly, 
by taking into account that the use of the regularity assumption \eqref{hppiuregdati} on the initial data 
would be essentially the same for the \pier{rigorous} procedure and the formal one.}

We differentiate \eqref{secondal} with respect to time and test the resulting equality by~$\dt\phil$.
At the same time, we test \eqref{primal} and \eqref{terzal} by $\dt\mul$ and~$\dt\Sl$, respectively.
Then, we sum up and integrate over~$(0,t)$.
The terms involving the product \juerg{\,$\dt\phil\,\dt\mul$\,} cancel each other, and we obtain
\Bsist
  && \alpha \intQt |\dt\mul|^2
  + \frac 12 \, \norma{A^\rho\mul(t)}^2
  \non
  \\
  && \quad {}
  + \frac \beta 2 \, \norma{\dt\phil(t)}^2
  + \iot \norma{B^\sigma\dt\phil(t)}^2
  + \intQt (\ful)'(\phil) |\dt\phil|^2
  \non
  \\
  && \quad {}
  + \pier{\intQt |\dt \Sl|^2}
  + \frac 12 \, \norma{C^\tau \juerg{\Sl}(t)}^2
  \non
  \\
  && = - \intQt f_2'(\phil) |\dt\phil|^2
  + \frac 12 \, \norma{A^\rho\muz}^2
  + \frac \beta 2 \, \norma{\dt\phil(0)}^2
  + \frac 12 \, \norma{C^\tau\Sz}^2
  \non
  \\
  && \quad {}
  + \intQt P(\phil) \bigl[ \juerg{(\Sl-\mul) \dt\mul - (\Sl-\mul) \dt \Sl} \bigr] \,.
  \label{perreg}
\Esist
All of the terms on the \lhs\ are nonnegative, and the first one on the \rhs\ is 
estimated by \pier{a constant proportional to $1/\beta$, due to \eqref{truestimate}}.
Moreover, the last integral, which we denote by $I$ for brevity, can be dealt with \pier{by using the Young inequality and~\eqref{truestimate}:
\begin{align}
   I \,&\leq \, c \intQt |(\Sl-\mul) \dt\mul| +
   c \intQt |(\Sl-\mul) \dt \Sl|
   \non
   \\
  \,&\leq\, 
   \frac \alpha 2 \intQt |\dt \mul|^2   
   + \frac 12 \intQt |\dt \Sl|^2   
   + c \tonde{ \frac1\alpha + 1} \iot \tonde{\norma{\Sl(s)}^2 +  \norma{\mul(s)}^2}
      \ds\,,     
  \non
  \\
  & \leq\, 
   \frac \alpha 2 \intQt |\dt \mul|^2   
   + \frac 12 \intQt |\dt \Sl|^2   
   + c_{\alpha}.
   \label{pier5}
\end{align}
It remains} to deal with the $H$-norm of \,$\dt\phi(0)$.
To this end, we observe that \eqref{eqsecondal} yields
\Beq
  \beta \,\dt\phil(0)
  = \muz - B^{2\sigma}\phiz - \ful(\phiz) - f_2(\phiz)\,,
  \non
\Eeq
whence (see \eqref{hppiuregdati} and \eqref{propyosidabis})
\Beq
  \norma{\dt\phil(0)}
  \leq \frac 1 \beta \bigl(
    \norma\muz 
    + \norma\phiz_{B,2\sigma}
    + \norma{f_1^\circ(\phiz)}
    + c ( \norma\phiz + 1)
  \bigr)
  \leq c_\beta \,.
  \non
\Eeq
Therefore, if we come back to~\eqref{perreg} \pier{and account for \eqref{pier5} and}  
the estimate \eqref{stimal} of the previous section, 
we see that we have proved that
\Beq
  \norma\mul_{\H1H\cap\L\infty{\VA\rho}}
  + \norma\phil_{\W{1,\infty}H\cap\H1{\VB\sigma}}
  + \norma S_{\H1H\cap\L\infty{\VC\tau}}
  \leq c_{\alpha,\beta} \,.
  \non
\Eeq
Since $(\mul,\phil,\Sl)$ converges to $(\mu,\phi,S)$ as shown in the previous section,
the above estimate implies many of the regularity properties stated in \accorpa{piuregmu}{piuregS}.
Indeed, by accounting for what we have already \pier{shown} in the first part,
we see that just the conditions 
$\mu\in\L2{\VA{2\rho}}$ \pier{and $S\in\L2{\VC{2\tau}}$ are missing. But
these properties immediately follow by comparison in the equations 
\eqref{prima} and~\eqref{terza}. This concludes the proof.}

\step
Proof of Theorem~\ref{Secondeq}

\gianni{In contrast to the proof of Theorem~\ref{Regularity},
we \juerg{here} use a completely rigorous argument since the details are not complicated.
However, a remark is necessary.
We recall that the assumption \eqref{embeddings} on the spaces $\VA\rho$ and $\VC\tau$
is not required in the statement, 
so that \juerg{uniqueness is neither ensured} for problem \Pbl\ nor for the approximating problem.
Hence, we have to be precise.
Namely, we fix any solution $\soluz$ that can be obtained by the procedure adopted in Section~\ref{EXISTENCE}
and prove both its further regularity and the existence of some $\xi$ satisfying the conditions of the statement
if \eqref{hpsecondeq} is fulfilled.
Thus, the interpolants of the discrete solution converge as $h$ tend to zero (along a suitable subsequence)
to~some solution to the approximating problem, which converges as $\lambda$ tends to zero (along a subsequence)
to~the solution we have chosen. 
So, coming back to the discrete problem \Pbld,
we multiply \eqref{secondad} by $\ful(\phinp)$.}
We notice that $\ful(\phinp)\in H$, due to \eqref{hpsecondeq}
with \,$v=\phinp$\, and \,$\psi=\ful$\,, 
\gianni{since $\phinp\in\VB{2\sigma}$, and \juerg{because} $\ful$ is monotone and \Lip\ continuous and vanishes at the origin}.
We obtain
\Bsist
  && \frac \beta h \, \bigl( \phinp-\phin , \ful(\phinp) \bigr)
  + \bigl( B^{2\sigma} \phinp , \ful(\phinp) \bigr)
  + \norma{\ful(\phinp)}^2
  \non
  \\
  && =  \bigl( \munp - f_2(\phinp) + L (\phin-\phinp) , \ful(\phinp) \bigr) \,.
  \label{pereqseconda}
\Esist
For the first term on the \lhs, we use the convexity in this way:
\Beq
  \bigl( \phinp-\phin , \ful(\phinp) \bigr)
  \geq \iO \Ful(\phinp) - \iO \Ful(\phin) \,.
  \non
\Eeq
The second term of \eqref{pereqseconda} is nonnegative by assumption~\eqref{hpsecondeq}.
Finally, the \rhs\ is estimated by owing to the Young inequality 
and to the linear growth of $f_2$ given by its \Lip\ continuity.
Namely, we have that
\Bsist
  && \bigl( \munp - f_2(\phinp) + L (\phin-\phinp) , \ful(\phinp) \bigr) 
  \non
  \\
  && \leq \frac 12 \, \norma{\ful(\phinp)}^2
  +  c \, (\norma\munp^2 + \norma\phinp^2 + \norma{\phin-\phinp}^2 + 1) \,.
  \non
\Esist
Therefore, combining with \eqref{pereqseconda}, rearranging, multiplying by~$h$, 
and summing up with respect to $n=0,\dots,N-1$, we deduce that
\Bsist
  && \beta \iO \Ful(\phi^N)
  + \frac 12 \, \somma n0{N-1} h\, \norma{\ful(\phinp)}^2
  \non
  \\
  && \leq \beta \iO \Ful(\phiz)
  + c \Bigl(
    \somma n0{N-1} h \,\norma\munp^2
    + \somma n0{N-1} h \,\norma\phinp^2
    + \somma n0{N-1} h \,\norma{\phin-\phinp}^2
  \Bigr).
  \non
\Esist
Since the first term on the \lhs\ is nonnegative,
the above inequality and {}{the second inequality in}~\eqref{propyosida}
imply the following estimate for the interpolants:
\Beq
  \norma{\ful(\overphih)}_{\L2H}^2
  \leq \beta\iO F_1(\phiz)
  + c \,\bigl(
    \norma\overmuh_{\L2H}^2
    + \norma\overphih_{\L2H}^2
    + \norma{\overphih-\underphih}_{\L2H}^2
  \bigr) .
  \non
\Eeq
By recalling \eqref{hpdati} for $\phiz$ and the estimates~\accorpa{basic}{primastimah},
we conclude that
\Beq
  \norma{\ful(\overphih)}_{\L2H} \leq c_{\alpha,\beta} \,.
  \non
\Eeq
Since $\overphih\to\phil$ strongly in $\L2H$ (see~\eqref{strongphih}) and $\ful$ is \Lip\ continuous, 
we infer that 
\Beq
  \norma{\ful(\phil)}_{\L2H} \leq c_{\alpha,\beta} \,.
  \non
\Eeq
Moreover, $\phil$~converges to $\phi$ strongly in $\L2H$.
Therefore, by using weak compactness and applying, e.g., \cite[Lemma~2.3, p.~38]{Barbu},
we conclude that
\Beq
  \ful(\phil) \to \xi
  \quad {}{\hbox{weakly in $\L2H$,\quad for some \,$\xi$\, with}} \quad
  \xi \in f_1(\phi)
  \quad \aeQ.
  \non
\Eeq
At this point, we can let $\lambda$ tend to zero in the integrated version \eqref{eqsecondal}
and deduce that
\Bsist
  && \beta \ioT \bigl( \dt\phi(t) , v(t) \bigr) \, dt
  + \ioT \bigl( B^\sigma\phi(t) , B^\sigma v(t) \bigr) \, dt
  + \ioT \bigl( \xi(t) + f_2(\phi(t) , v(t) \bigr) \, dt
  \non
  \\
  && = \ioT \bigl( \mul(t) , v(t) \bigr) \, dt
  \quad \hbox{for every $v\in\L2{\VB\sigma}$}.
  \non
\Esist
This variational equation is equivalent to 
\Beq
  \beta \,\dt\phi + B^{2\sigma}\phi + \xi + f_2(\phi) = \mu
  \quad \hbox{\aet\ in the sense of $\VB{-\sigma}$},
  \non
\Eeq
and this implies both \eqref{regseconda} and \eqref{eqseconda}.
\gianni{In order to prove the last sentence, it suffices to recall that \juerg{the}
\pier{embedding properties \eqref{embeddings} ensure} uniqueness for the solution $\soluz$.
Hence the uniqueness of $\xi$ simply follows by comparison in~\eqref{eqseconda}.}

\step
Proof of Corollary~\ref{Fullregularity}

\gianni{The assumptions of the statement guarantee that 
the solution $\soluz$ is unique and that there exists a unique $\xi$
satisfying the properties stated in Theorem~\ref{Secondeq}.
In particular, by the above proofs,
the (unique) solution $\soluzl$ to the approximating problem and the corresponding $\ful(\phil)$
converge to $\soluz$ and to~$\xi$, respectively,
in the proper topologies.
Moreover, as $\ful$ satisfies the same assumptions \juerg{as those} we have postulated for~$f_1$,
we can apply Theorem~\ref{Secondeq} to the approximating problem.
Hence, $\phil$ belongs to $\L2{\VB{2\sigma}}$, and \eqref{secondal} can be replaced by the equation
\Beq
  B^{2\sigma} \phil(t) + \ful(\phil(t)) 
  = \mul(t) - \beta \, \dt\phil(t) - f_2(\phil(t))
  \quad \aat.
  \label{strongsecondal}
\Eeq
For a while, we argue for a fixed~$t$ (\aet).
We multiply \eqref{strongsecondal} by $\ful(\phi(t))\in H$.
Since $\phil(t)\in\VB{2\sigma}$ and $\ful$ is monotone, \Lip\ continuous, and vanishes at the origin,
we can apply \eqref{hpsecondeq} and have that
\Beq
  \bigl( B^{2\sigma} \phil(t) , \ful(\phil(t)) \bigr) \geq 0 \,.
  \non
\Eeq
Therefore, we obtain the inequality
\Beq
  \norma{\ful(\phil(t))}
  \leq \norma{\mul(t) - \beta \, \dt\phil(t) - f_2(\phil(t))} \,.
  \non
\Eeq
Since $\mul-\beta\,\dt\phil-f_2(\phil)$ is bounded in $\L\infty H$ by \eqref{truestimate},
the same is true for $\ful(\phil)$.
By comparison in~\eqref{strongsecondal}, we deduce that $B^{2\sigma}\phil$ is bounded as well.
Hence, it immediately follows that $B^{2\sigma}\phi\in\L\infty H$,
whence also $\xi\in\L\infty H$ by comparison in~\eqref{eqseconda}.
}


\section*{Acknowledgments}
\pier{This research was supported by the Italian Ministry of Education, 
University and Research~(MIUR): Dipartimenti di Eccellenza Program (2018--2022) 
-- Dept.~of Mathematics ``F.~Casorati'', University of Pavia. 
In addition, PC and CG gratefully acknowledge some other 
financial support from the GNAMPA (Gruppo Nazionale per l'Analisi Matematica, 
la Probabilit\`a e le loro Applicazioni) of INdAM (Isti\-tuto 
Nazionale di Alta Matematica) and the IMATI -- C.N.R. Pavia.}


\vspace{3truemm}

\Begin{thebibliography}{10}

\bibitem{Barbu}
V. Barbu,
``Nonlinear Differential Equations of Monotone Type in Banach Spaces'',
Springer,
London, New York, 2010.

\bibitem{Brezis}
H. Brezis,
``Op\'erateurs maximaux monotones et semi-groupes de contractions
dans les espaces de Hilbert'',
North-Holland Math. Stud.
{\bf 5},
North-Holland,
Amsterdam,
1973.

\bibitem{BLM}
N. Bellomo, N.~K. Li, P.~K. Maini,
On the foundations of cancer modelling: selected topics, speculations, and perspectives,
\textit{Math. Models Methods Appl. Sci.} \textbf{18} (2008), 593--646.

\bibitem{BCG}
S. Bosia, M. Conti, M. Grasselli,
On the Cahn--Hilliard--Brinkman system,
\textit{Commun. Math. Sci.} \textbf{13} (2015), 1541--1567.

\bibitem{Cahn}
J.~W. Cahn, J.~E. Hilliard,
Free energy of a nonuniform system. I. Interfacial free energy,
{\it J. Chem. Phys.} {\bf 28} (1958), 258--267.

\bibitem{CRW}
C. Cavaterra, E. Rocca, H. Wu,
Long-time dynamics and optimal control of a diffuse interface model for tumor growth,
\textit{Appl. Math. Optim.} (2019), Online First 15 March 2019, https://doi.org/10.1007/s00245-019-09562-5

\bibitem{CWSL}
Y. Chen, S.~M. Wise, V.~B. Shenoy, J.~S. Lowengrub,
A stable scheme for a nonlinear multiphase tumor growth model with an elastic membrane,
\textit{Int. J. Numer. Methods Biomed. Eng.} \textbf{30} (2014), 726--754.

\pier{\bibitem{CG1} P. Colli, G. Gilardi, 
Well-posedness, regularity and asymptotic analyses 
for a fractional phase field system,
{\em Asymptot. Anal.}, to appear (see also preprint 
arXiv:1806.04625 [math.AP] (2018), pp.~1--34).}

\bibitem{CGH15}
P. Colli, G. Gilardi, D. Hilhorst,
On a Cahn--Hilliard type phase field system related to tumor growth,
\textit{Discrete Contin. Dyn. Syst.} \textbf{35} (2015), 2423--2442.

\pier{\bibitem{CGMR}
P. Colli, G. Gilardi, G. Marinoschi, E. Rocca, 
Sliding mode control for a phase field system related to tumor growth.
Appl. Math. Optim. \textbf{79} (2019), 647--670.}

\bibitem{CGRS1}
P. Colli, G. Gilardi, E. Rocca, J. Sprekels,
Vanishing viscosities and error estimate for a Cahn--Hilliard type phase field system related to tumor growth,
\textit{Nonlinear Anal. Real World Appl.} \textbf{26} (2015), 93--108.

\bibitem{CGRS2}
P. Colli, G. Gilardi, E. Rocca, J. Sprekels, 
Asymptotic analyses and error estimates for a Cahn--Hilliard type phase field system modelling tumor growth,
\textit{Discrete Contin. Dyn. Syst. Ser. S.} {\bf  10} (2017), 37--54. 

\bibitem{CGRS3}
P. Colli, G. Gilardi, E. Rocca, J. Sprekels,
Optimal distributed control of a diffuse interface model of tumor growth,
\textit{Nonlinearity} \textbf{30} (2017), 2518--2546.

\bibitem{CGS18}
P. Colli, G. Gilardi, J. Sprekels,
Well-posedness and regularity for a generalized fractional Cahn--Hilliard system,
{\em Atti Accad. Naz. Lincei Rend. Lincei Mat. Appl.}, to appear 
(see also preprint arXiv:1804.11290 [math.AP] (2018), pp. 1--36).

\bibitem{CGS19}
P. Colli, G. Gilardi, J. Sprekels,
Optimal distributed control of a generalized fractional Cahn--Hilliard system,
{\em Appl. Math. Optim.}, Online First 15 November 2018, https://doi.org/10.1007/s00245-018-9540-7.

\bibitem{CGS22}
P. Colli, G. Gilardi, J. Sprekels,
Longtime behavior for a generalized Cahn--Hilliard system with fractional operators,
preprint arXiv:1904.00931 [math.AP] (2019), pp.~1--18.

\bibitem{ConGio}
M. Conti, A. Giorgini,
The three-dimensional Cahn--Hilliard--Brinkman system with unmatched viscosities,
\pier{preprint hal-01559179 (2018), pp.~1--34.}

\pier{\bibitem{CLLW}
V. Cristini, X. Li, J.~S. Lowengrub, S.~M. Wise,
Nonlinear simulations of solid tumor growth using a mixture model: invasion and branching,
{\it J. Math. Biol.} {\bf 58} (2009), 723--763.}

\bibitem{CL2010}
V. Cristini, J.~S. Lowengrub,
``Multiscale Modeling of Cancer: an Integrated Experimental and Mathematical Modeling Approach'',
Cambridge Univ. Press, Cambridge, 2010.

\bibitem{DFRSS}
M. Dai, E. Feireisl, E. Rocca, G. Schimperna, M. Schonbek,
Analysis of a diffuse interface model for multi-species tumor growth,
\textit{Nonlinearity} \textbf{30} (2017), \pier{1639--1658.}

\bibitem{DG}
F. Della Porta, M. Grasselli,
On the nonlocal Cahn--Hilliard--Brinkman and Cahn--Hilliard--Hele--Shaw systems,
\textit{Commun. Pure Appl. Anal.} \textbf{15} (2016), 299--317,
Erratum: \textit{Commun. Pure Appl. Anal.} \textbf{16} (2017), 369--372.
%

\bibitem{EGAR}
M. Ebenbeck, H. Garcke,
Analysis of a Cahn--Hilliard--Brinkman model for tumour growth with chemotaxis,
{\it J. Differential Equations} \textbf{266} (2019), 5998--6036.

\bibitem{FBG2006}
A. Fasano, A. Bertuzzi, A. Gandolfi,
Mathematical modeling of tumour growth and treatment,
\pier{in ``Complex Systems in Biomedicine'', 
A.~Quarteroni, L.~Formaggia, A.~Veneziani~(ed.),
Springer, Milan,  2006, pp.~71--108.}

\bibitem{FW2012}
X. Feng, S.~M. Wise,
Analysis of a Darcy--Cahn--Hilliard diffuse interface model for the Hele--Shaw flow and its fully discrete finite element approximation,
\textit{SIAM J. Numer. Anal.} \textbf{50} (2012), 1320--1343.

\bibitem{Fri2007}
A. Friedman,
Mathematical analysis and challenges arising from models of tumor growth,
\textit{Math. Models Methods Appl. Sci.} \textbf{17} (2007), 1751--1772.

\bibitem{Lowen10}
H.~B. Frieboes, F. Jin, Y.~L. Chuang, S.~M. Wise, J.~S. Lowengrub, V. Cristini,
Three-dimensional multispecies nonlinear tumor growth - II: tumor invasion and angiogenesis,
{\it J. Theoret. Biol.} {\bf 264} (2010), 1254--1278.

\bibitem{FGR}
S. Frigeri, M. Grasselli, E. Rocca,
On a diffuse interface model of tumor growth,
\textit{European J. Appl. Math.} \textbf{26} (2015), 215--243.

\pier{\bibitem{FLR}
S. Frigeri, K.~F. Lam, E. Rocca, 
On a diffuse interface model for tumour growth with non-local 
interactions and degenerate mobilities, 
in ``Solvability, regularity, and optimal control of boundary 
value problems for PDEs'', 
P.~Colli, A.~Favini, E.~Rocca, G.~Schimperna, J.~Sprekels~(ed.), 
Springer INdAM Series~{\bf 22}, Springer, Cham, 2017, pp.~217--254.}

\bibitem{FLRS}
S. Frigeri, K.~F. Lam, E. Rocca, G. Schimperna,
On a multi-species Cahn--Hilliard--Darcy tumor growth model with singular potentials,
\textit{Commun. Math Sci.} \textbf{16} (2018), 821--856.

\bibitem{GL2016}
H. Garcke, K.~F. Lam,
Global weak solutions and asymptotic limits of 
a Cahn--Hilliard--Darcy system modelling tumour growth,
\textit{AIMS Mathematics} \textbf{1} (2016), 318--360.

\pier{%
\bibitem{GL2017-1}
H. Garcke,  K.~F. Lam,
Analysis of a Cahn--Hilliard system with non--zero Dirichlet 
conditions modeling tumor growth with chemotaxis,
{\it Discrete Contin. Dyn. Syst.} {\bf 37} (2017), 4277--4308.
\bibitem{GL2017-2}
H. Garcke, K.~F. Lam,
Well-posedness of a Cahn--Hilliard system modelling tumour
growth with chemotaxis and active transport,
{\it European. J. Appl. Math.} {\bf 28} (2017), 284--316.}

\bibitem{GL2018}
H. Garcke, K.~F. Lam,
On a Cahn--Hilliard--Darcy system for tumour growth with solution dependent source terms,
in ``Trends on Applications of Mathematics to Mechanics'', 
E.~Rocca, U.~Stefanelli, L.~Truskinovski, A.~Visintin~(ed.), 
Springer INdAM Series~{\bf 27}, Springer, Cham, 2018, pp.~243--264.

\bibitem{GLNS}
H. Garcke, K.~F. Lam, R. N\"urnberg, E. Sitka,
A multiphase Cahn--Hilliard--Darcy model for tumour growth with necrosis,
\textit{Math. Models Methods Appl. Sci.} \textbf{28} (2018), 525--577.

\bibitem{GLR}
H. Garcke, K.~F. Lam, E. Rocca,
Optimal control of treatment time in a diffuse interface model for tumour growth,
\textit{Appl. Math. Optim.} {\bf 78} (2018), 495--544.

\bibitem{GLSS}
H. Garcke, K.~F. Lam, E. Sitka, V. Styles,
A Cahn--Hilliard--Darcy model for tumour growth with chemotaxis and active transport,
\textit{Math. Models Methods Appl. Sci.} \textbf{26} (2016), 1095--1148.

\bibitem{GioGrWu}
A. Giorgini, M. Grasselli, H. Wu,
The Cahn--Hilliard--Hele--Shaw system with singular potential,
\textit{Ann. Inst. H. Poincar\'e Anal. Non Lin\'eaire} 
{\bf 35} (2018), 1079--1118.

\pier{
\bibitem{HDPZO}
A. Hawkins-Daarud, S. Prudhomme, K.~G. van der Zee, J.~T. Oden,
Bayesian calibration, validation, and uncertainty quantification of diffuse 
interface models of tumor growth,
{\it J. Math. Biol.} {\bf 67} (2013), 1457--1485.}

\bibitem{HZO12}
A. Hawkins-Daarud, K.~G. van der Zee, J.~T. Oden,
Numerical simulation of a thermodynamically consistent four-species tumor growth model,
\textit{Int. J. Numer. Meth. Biomed. Engrg.} \textbf{28} (2012), 3--24.

\pier{\bibitem{HKNZ}
D. Hilhorst, J. Kampmann, T.~N. Nguyen, K.~G. van der Zee, Formal asymptotic
limit of a diffuse-interface tumor-growth model, 
{\it Math. Models Methods Appl. Sci.} {\bf 25} (2015), 1011--1043.}

\bibitem{JWZ}
J. Jiang, H. Wu, S. Zheng,
Well-posedness and long-time behavior of a non-autonomous
Cahn--Hilliard--Darcy system with mass source modeling tumor growth,
{\it J. Differential Equations} {\bf 259} (2015), 3032--3077.

\bibitem{Lions}
J.-L. Lions, ``Quelques M\'ethodes de R\'esolution des Probl\`emes aux Limites non Lin\'eaires'', 
Dunod; Gauthier-Villars, Paris, 1969. 

\bibitem{LTZ}
J.~S. Lowengrub, E.~S. Titi, K. Zhao,
Analysis of a mixture model of tumor growth,
{\it European J. Appl. Math.} \textbf{24} (2013), 691--734.
%
%

\pier{\bibitem{MRS}
A. Miranville, E. Rocca, G. Schimperna,
On the long time behavior of a tumor growth model,
{\it J. Differential Equations\/} {\bf 267} (2019) 2616--2642.
\bibitem{OHP}
J.~T. Oden, A. Hawkins, S. Prudhomme,
General diffuse-interface theories and an approach to predictive tumor growth modeling,
{\it Math. Models Methods Appl. Sci.} {\bf 20} (2010), 477--517. 
\bibitem{Sig}
A. Signori,
Optimal distributed control of an extended model of
tumor growth with logarithmic potential,
{\em Appl. Math. Optim.}, Online First 30 October 2018, 
https://doi.org/10.1007/s00245-018-9538-1.
\bibitem{S_a}
A. Signori,
Optimal treatment for a phase field system of Cahn--Hilliard 
type modeling tumor growth by asymptotic scheme,
preprint arXiv:1902.01079 [math.AP] (2019), pp.~1--28.
\bibitem{S_b}
A. Signori,
Vanishing parameter for an optimal control problem modeling tumor growth,
preprint arXiv:1903.04930 [math.AP] (2019), pp.~1--22.}

\juergen{
\bibitem{SW}
J. Sprekels, H. Wu, Optimal distributed control of a Cahn--Hilliard--Darcy system with mass sources,
{\em Appl. Math. Optim.}, Online First 24 January 2019,
https://doi.org/10.1007/s00245-019-09555-4.
}

\bibitem{WW2012}
X.-M. Wang, H. Wu,
Long-time behavior for the Hele--Shaw--Cahn--Hilliard system,
\textit{Asymptot. Anal.} \textbf{78} (2012), 217--245.

\bibitem{WZ2013}
X.-M. Wang, Z.-F. Zhang,
Well-posedness of the Hele--Shaw--Cahn--Hilliard system,
\textit{Ann. Inst. H. Poincar\'e Anal. Non Lin\'eaire} \textbf{30} (2013), 367--384.

\bibitem{Wise2011}
S.~M. Wise, J.~S. Lowengrub, V. Cristini,
An adaptive multigrid algorithm for simulating solid tumor growth using mixture models,
\textit{Math. Comput. Modelling} \textbf{53} (2011), 1--20.

\bibitem{WLFC}
S.~M. Wise, J.~S. Lowengrub, H.~B. Frieboes, V. Cristini,
Three-dimensional multispecies nonlinear tumor growth - I: model and numerical method,
{\it J. Theoret. Biol.} {\bf 253} (2008), 524--543.

\pier{\bibitem{WZZ}
X. Wu, G.~J. van Zwieten, K.~G. van der Zee, Stabilized second-order splitting
schemes for Cahn--Hilliard~models with applications to 
diffuse-interface tumor-growth models, 
{\it Int. J. Numer. Methods Biomed. Eng.} {\bf 30} (2014), 180--203.}

\End{thebibliography}

\End{document}
